\documentclass[11pt]{amsart}

\usepackage{amsmath}
\usepackage{amssymb}
\usepackage{amscd}
\usepackage{amscd}
\usepackage{array}
\usepackage[all]{xy}

\def\c{{\rm C}}
\def\N{{\mathbb N}}
\def\C{{\mathbb C}}

\def\z{{\mathbb Z}}

\begin{document}

\newtheorem{theorem}{Theorem}[section]
\newtheorem{lemma}[theorem]{Lemma}
\newtheorem{corollary}[theorem]{Corollary}
\newtheorem{definition}[theorem]{Definition}
\newtheorem{proposition}[theorem]{Proposition}
\newtheorem{theo}[theorem]{Theorem}
\newtheorem{prop}[theorem]{Proposition}
\newtheorem{coro}[theorem]{Corollary}
\newtheorem{lemm}[theorem]{Lemma}
\newtheorem{defi}[theorem]{Definition}
\newtheorem{rem}[theorem]{Remark}
\newtheorem{remark}[theorem]{Remark}
\newtheorem{ex}[theorem]{Example}
\newtheorem{example}[theorem]{Example}
\newtheorem{question}[theorem]{Question}

\title{hopf images and inner faithful representations}

\author{Teodor Banica}
\address{D\'epartement de Mathématiques, Universit\'e de Cergy-Pontoise, 2 avenue Chauvin, 95302 Cergy-Pontoise, France} 
\email{Teodor.Banica@u-cergy.fr}

\author{Julien Bichon}
\address{Laboratoire de Math\'ematiques, Universit\'e Blaise Pascal, Clermont-Ferrand II, 
Campus des C\'ezeaux, 63177 Aubi\`ere Cedex, France}
\email{Julien.Bichon@math.univ-bpclermont.fr}

\thanks{Work supported by the ANR project ''Galoisint'' BLAN07-3\_183390}

\subjclass[2000]{16W30}

\begin{abstract}
We develop a general theory of Hopf image of a Hopf algebra representation, with  the associated concept of inner faithful representation, modelled on the notion of faithful representation of a discrete group. We study several examples, including group algebras, enveloping algebras of Lie algebras, pointed Hopf algebras, function algebras, twistings and cotwistings, and we present a Tannaka duality formulation of the notion of Hopf image.
\end{abstract}

\maketitle

\section{Introduction}
The aim of this paper is to provide an axiomatization and systematic study of the concept of Hopf image of a Hopf algebra representation, as well as the related concept of inner faithful representation. These notions appeared, under various degrees of generality, in a number of independent investigations: vertex models and related quantum groups \cite{ba98}, \cite{bm}, \cite{bn}, locally compact quantum groups and their outer actions \cite{kv}, \cite{va}. 
These were used extensively in the recent paper \cite{bbs}
in order to study the quantum symmetries of Hadamard matrices and of the corresponding subfactors.

The leading idea is that we want to translate the notion of faithful representation of a discrete
group at a Hopf algebra level. Let $\Gamma$ be a group, let 
$A$ be a $k$-algebra ($k$ is a field) and consider a representation
$$\pi : k[\Gamma] \longrightarrow A$$
There are two possible notions of faithfulness for $\pi$:
\begin{enumerate}
 \item faithfulness of $\pi$ as an algebra map, in which case we simply say that $\pi$ is faithful;
 \item faithfulness of the induced group morphism $\pi_{|\Gamma} :  \Gamma \longrightarrow A^\times$.
\end{enumerate}
It is clear that the first notion is much more restrictive than the second one, and 
choice of one of these notions as the good one for faithfulness
depends on whether one is interested by the algebra $k[\Gamma]$
or rather by the group $\Gamma$ itself.

It is not difficult so see that $\pi_{|\Gamma}$ is faithful if and only if
${\rm Ker}(\pi)\subset k[\Gamma]$ does not contain any non-zero Hopf ideal. 
This simple observation leads to a notion of inner faithful
representation for arbitrary Hopf algebras: if $H$ a Hopf algebra,  we say that
a representation
 $$\pi : H \longrightarrow A$$
 is \textsl{inner faithful} is ${\rm Ker}(\pi)$ does not contain any non-zero Hopf ideal.
 Of course $H$ is viewed as  the group algebra
 of a discrete quantum group, and faithfulness refers to this discrete quantum group.
 
 If the representation $\pi$ is not inner faithful, there is however a minimal
 Hopf algebra $H_\pi$ that factorizes $\pi$, that we call the \textsl{Hopf image} of $\pi$
 (for $H=k[\pi]$, we have $k[\Gamma]\pi= k[\pi(\Gamma)]$). 
 The exact universal property of the Hopf image is stated in Section 2.
The Hopf image measures how much $\pi$ fails to be faithful in the discrete
quantum group sense. In the situation of \cite{bbs}, the Hopf image
of the representation of the quantum permutation algebra
associated to a complex Hadamard matrix  measures the complexity of the quantum invariants of the Hadamard matrix.

 A natural concept arising from these considerations is the notion
 of $\textsl{inner linear}$ Hopf algebra: we say that a Hopf algebra is inner
 linear if it admits a finite-dimensional inner faithful 
representation. Therefore the problem of inner linearity
for Hopf algebras is a generalization of the celebrated
linearity problem for discrete groups.
More generally, we believe that the study of Hopf images leads 
to interesting new questions and problems in Hopf algebra theory, as well as in group
theory through the study of group duals.

The paper is devoted to a general study of the notion of Hopf image, culminating in a Tannakian
formulation, together with the study of several key examples.

The precise contents of the paper is as follows. 
In Section 2 we give a precise formulation of the concept of Hopf image,
we prove the existence of the Hopf image 
and study  some of its basic properties.
In Section 3, we get back to the discrete
group example discussed in the introduction and also discuss
the example of enveloping algebras of Lie algebras. 
We also introduce the notion of inner linear Hopf algebra.
In Section 4  
we consider the case of pointed Hopf algebras (thus including quantized enveloping algebras), for 
which we give a criterion of inner faithfulness for a representation using skew-primitives. As a simple illustration, we
examine the vector  $2$-dimensional representation
of $U_q(\mathfrak{sl}_2)$, which is shown to be inner faithful if and only if $q$ is not a root of unity.
Section 5 is devoted to function algebras on algebraic and compact groups.
The construction of inner faithful representations is related to
the existence of topological generators.
These simple examples already show that for cosemisimple Hopf algebras, inner
faithulness is not directly related to injectivity on characters (in contrast with the group algebra situation).
 In Section 6 we consider the problem of constructing 
 inner faithful representations for 
Hopf algebras with twisted coproduct or twisted product.
We show that the cotwist ($2$-cocycle deformation)
of an inner linear Hopf algebra for which $S^2$ is an inner automorphism
is still inner linear if the 2-cocycle is
induced by a finite-dimensional quotient.
In particular  the multiparametric 2-cocycle deformations (cotwists) of 
$\mathcal O({\rm GL}_n(\C))$ at roots of unity are inner linear.
Section 7 is devoted to the study the tensor and free products of Hopf images.
The question of knowing whether a tensor product 
of inner faithful representations is still inner faithful leads to the notion of projectively inner faithful representation.
Section 8 is devoted to Tannaka duality type results, which
provide so far the best inner faithfulness criteria for a representation of a compact Hopf algebra
(i.e. a Hopf algebra associated to a compact quantum group).
These are used in Section 9 to describe examples of inner faithful
representations of  compact Hopf algebras having all their
simple comodules of dimension smaller than $2$. 

\medskip

\noindent
\textbf{Notations and conventions.} We work in general over a fixed field $k$.
We assume that the reader has  some familiarity with Hopf algebras, for which
the textbooks \cite{mo} or \cite{sw} are convenient. Our terminology and notation are the standard ones: in particular, for a Hopf algebra, $\Delta$, $\varepsilon$ and $S$ denote the comultiplication, counit and antipode, respectively.  
 
\medskip

\noindent
\textbf{Ackowledgements.} We wish to thank Rupert Yu for helpful discussions
on algebraic groups.

\section{Construction of the Hopf image and basic properties}

In this section we give the precise formulation of the concept of Hopf image,
prove its existence 
and study of some of its basic properties. The case of Hopf $*$-algebras is also considered
at the end of the section.

\subsection{Construction of the Hopf image}
First, let us give a precise formulation of the notion of Hopf image of a representation. 
Let $H$ be a Hopf algebra over a field $k$. As usual, a representation of
$H$ on an algebra $A$ is an algebra morphism $\pi : H \longrightarrow A$.

Let us say that a \textbf{factorization of $\pi$}
is a triple $(L, q, \varphi)$  where $L$ is a Hopf algebra,
$q : H \longrightarrow L$ is a surjective Hopf algebra map
and $\varphi : L \longrightarrow A$ is a representation, with 
the decomposition $\pi = \varphi \circ q$.
We define in a straightforward manner the category of factorizations
of $\pi$, and the \textbf{Hopf image of $\pi$} is defined to be the
final object in this category (hence we can also say that this is a minimal factorization).

\begin{theorem}
Let $\pi : H \longrightarrow A$ be a representation
of a Hopf algebra $H$ on an algebra $A$. Then $\pi$
has a Hopf image: 
there exists a Hopf algebra $H_\pi$ together
with a surjective Hopf algebra map $p : H \longrightarrow H_\pi $ and 
 a representation $\tilde{\pi} :  H_\pi  \longrightarrow A$, such that
$\pi = \tilde{\pi} \circ p$, and such that if $(L,q, \varphi)$ is another 
factorization of $\pi$, there exists a unique
Hopf algebra map $f : L \longrightarrow H_\pi$ such that
$f \circ q = p$ and $\tilde{\pi} \circ f = \varphi$.  
$$\xymatrix{
H \ar[rr]^{\pi} \ar[dr]_{p} \ar@/_/[ddr]_q & & A  \\
& H_\pi \ar[ur]_{\tilde{\pi}} & \\
& L \ar@{-->}[u]_f \ar@/_/[uur]_{\varphi}& 
}$$
\end{theorem}

\begin{proof}  Let $I_\pi$ be sum of all the Hopf ideals contained in 
in ${\rm Ker}(\pi)$. It is clear that $I_\pi$ is a Hopf ideal
and is  the largest Hopf ideal contained in ${\rm Ker}(\pi)$.
 Let $H_\pi = H/I_\pi$ and 
 let $p : H \longrightarrow H_\pi$ be the canonical surjection. 
 Since $I_\pi \subset {\rm Ker}(\pi)$, there exists
 a unique algebra map $\tilde{\pi} : H_\pi \longrightarrow A$ such that $\tilde{\pi} \circ p = 
 \pi$.

 Let $(L,q, \varphi)$ be a factorization of $\pi$.
  Then ${\rm Ker}(q)$ is a Hopf ideal contained in 
  ${\rm Ker}(\pi)$ and hence ${\rm Ker}(q) \subset I_\pi$. Hence
 there exists a unique Hopf algebra map $f : L \longrightarrow H_\pi$, $q(x) \longmapsto p(x)$, satisfying
 $f \circ q=p$ and $\varphi \circ f = \tilde{\pi}$, and thus $H_\pi$ has the required universal property.
\end{proof}

The above proof uses the Hopf ideal $I_\pi$, the largest Hopf ideal
contained in ${\rm Ker}(\pi)$. This Hopf ideal is constucted in a very abstract manner, and it is
useful for several purposes to have a more concrete description of $I_\pi$, that we give now.

Let $F$ be the free monoid generated by the set $\mathbb N$, with its
generators denoted $\alpha_0, \alpha_1 ,\ldots$ and its unit element denoted $1$.
We denote by $\tau : F \longrightarrow F$ the unique monoid antimorphism
such that $\tau(\alpha_i) = \alpha_{i+1}$.

To any element $g \in F$, we associate an algebra $A^g$, defined inductively on the lenght
of $g$ as follows.
We put $A^1=k$, $A^{\alpha_k} = A$ if $k$ is even and $A^{\alpha_k} = A^{\rm op}$
if $k$ is odd. Now for $g,h \in F$ with $l(g) > 1$ and $l(h)>1$, we put
$A^{gh}= A^g \otimes A^h$.

Now we associate an algebra morphism $\pi^g : H \longrightarrow A^g$ to any
$g \in F$, again by induction on the length of $g$.  We put $\pi^1 = \varepsilon$, $\pi^{\alpha_k} = \pi \circ S^k$, and 
for $g,h \in F$ with $l(g)>1$ and $l(h)>1$, we put $\pi^{gh} = (\pi^g \otimes \pi^h) \circ \Delta$.

\begin{proposition} Let $\pi : H \longrightarrow A$ be a representation
and let $I_\pi$ be the largest Hopf ideal contained in ${\rm Ker}(\pi)$. 
We have
$$I_\pi = \bigcap_{g \in F} {\rm Ker}(\pi^g)\subset H$$
\end{proposition}

The proof of Proposition 2.2 uses several lemmas. We put $J_\pi = \bigcap_{g \in F} {\rm Ker}(\pi^g)\subset H$.
By construction $J_\pi$ is an ideal in $H$, and
we wish to prove now that $J_\pi$ is a Hopf ideal. The following result ensures that it is a coideal.

\begin{lemma}
 Let $C_\pi$ be the linear subspace in $H^*$
 generated by the elements 
 $$\{ \psi \circ \pi^g, \ g \in F, \ \psi \in A_g^*\} \subset H^*$$ Then
 $C_\pi$ is a subalgebra of $H^*$, and $J_\pi = C_\pi^\perp$.
In particular $J_\pi$ is a coideal in $H$. 
\end{lemma}

\begin{proof}
 For $\psi,\phi \in H^*$ and $g,h \in F$, we have
 $$(\psi \circ \pi^g) \cdot (\phi \circ \pi^h)=
 ((\psi \circ \pi^g) \otimes (\phi \circ \pi^h)) \circ \Delta
 = (\psi \otimes \phi) \circ \pi^{gh}$$
 and since $\varepsilon \in C_\pi$, we conclude that $C_\pi$ is a subalgebra of $H^*$.
 It is clear that $C_\pi^\perp = I_\pi$, and we conclude 
 that $J_\pi$ is a coideal by Proposition 1.4.6 in \cite{sw}.
\end{proof}

\begin{lemma}
 For all $g \in F$, there exists a linear isomorphism
 $R_g : A^{\tau(g)} \longrightarrow A^g$ such that
 $\pi^g \circ S = R_g \circ \pi^{\tau(g)}$. In particular
 $S(J_\pi) \subset J_\pi$, and $J_\pi$ is a Hopf ideal in $H$.
\end{lemma}

\begin{proof}
We prove the result by induction on $l(g)$. For $l(g) =0$ or $l(g)=1$, this
follows immediately from the definitions. So assume that $l(g)>1$, so that
$g= hh'$, with $l(g)>l(h)\geq 1$ and $l(g)>l(h') \geq 1$. For $x$ in $H$, we have

\begin{align*}
 \pi^g \circ S(x) & = \pi^{hh'} \circ S(x)  = \pi^h(S(x_2)) \otimes \pi^{h'}(S(x_1)) = 
 ((R_h \circ \pi^{\tau(h)})(x_2)) \otimes (R_{h'} \circ \pi^{\tau(h')}(x_1)) \\
 & = R_{hh'} \circ (\pi^{\tau(h')} \otimes \pi^{\tau(h)}) \circ \Delta(x)
 =R_{g} \circ \pi^{\tau(g)}(x)
\end{align*}
The last assertion is now clear.
\end{proof}

The following lemma finishes the proof of Proposition 2.2.

\begin{lemma}
Let  $\pi : H \longrightarrow A$ be an algebra map, let
$q : H \longrightarrow L$ be a Hopf algebra map, and let
$\varphi : L \longrightarrow A$ be an algebra map with
$\pi = \varphi \circ q$. Then ${\rm Ker}(q) \subset J_\pi$, and hence 
any Hopf ideal contained in  ${\rm Ker}(\pi)$ is contained in $J_\pi$.
\end{lemma}

\begin{proof} Let us show 
 that $$\pi^g = \varphi^g \circ q, \ \forall g \in F$$
 with the same notation for $\varphi$ as the one for $\pi$.
 We show this by induction on $l(g)$.
 This is true if $l(g)=0$ since $q$ is a coalgebra map and
 if $l(g)=1$, we have
 $\pi^{\alpha_k} = \pi \circ S^k = \varphi \circ q \circ  S^k
 = \varphi \circ S^k \circ q = \varphi^{\alpha_k} \circ q$. 
 Assume now that $l(g) > 1$, so that $g=hh'$ with $l(g)>l(h)\geq 1$ and $l(g)>l(h') \geq 1$, and that the result is proved for elements of lenght $<l(g)$.
 Then 
 $$\pi^{hh'} = (\pi^h \otimes \pi^{h'}) \circ \Delta = 
 (\varphi^h \otimes \varphi^{h'}) \circ (q \otimes q) \circ \Delta
 = (\varphi^h \otimes \varphi^{h'}) \circ \Delta \circ q =
 \varphi^{hh'} \circ q$$
 and this proves our assertion by induction. Hence we have ${\rm Ker}(q) \subset J_\pi$.
 Any Hopf ideal is the kernel of a Hopf algebra map, and hence the last assertion follows.
\end{proof}

\begin{remark}
 {\rm 
 The notion of Hopf image considered here is  
 in general different from the one given in \cite{ad}, Definition 1.2.0, which
 refers to the category of Hopf algebras. Our
 definition of Hopf image uses the largest category of algebras.}   
\end{remark}

\subsection{Inner faithful representations}
The following definition was already given in the introduction.

\begin{defi}
Let $\pi : H \longrightarrow A$ be a representation
of a Hopf algebra $H$ on an algebra $A$. We say that 
$\pi$ is \textbf{inner faithful} if ${\rm Ker}(\pi)$ does not contain any non-zero Hopf ideal.
\end{defi}

We have several equivalent formulations for
inner faithfulness.

\begin{proposition}
 Let $\pi : H \longrightarrow A$ be a representation
of a Hopf algebra $H$ on an algebra $A$.
The following assertions are equivalent.
\begin{enumerate}
\item $\pi$ is inner faithful. 
\item
$\bigcap_{g \in F} {\rm Ker}(\pi^g) = (0)$
\item The Hopf algebra morphism $p : H \longrightarrow H_\pi$ in Theorem 2.1 is an isomorphism. 
\item If $(L,q,\varphi)$ is any factorization of $\pi$, then $q$ is an isomorphism.
\end{enumerate}
\end{proposition}

The equivalence between these assertions follows from the previous considerations and Proposition 2.2. 
Under a special assumption, we also have another
equivalent criterion for inner faithfulness.

\begin{proposition}
 Let $\pi : H \longrightarrow A$ be a representation
of a Hopf algebra $H$ on an algebra $A$ such that 
${\rm Ker}(\pi)$ is $S$-stable: $S({\rm Ker}(\pi)) \subset {\rm Ker}(\pi)$.
Then $\pi$ is inner faithful if and only if $(0)$ is the only bi-ideal contained
in ${\rm Ker}(\pi)$.
\end{proposition}

\begin{proof} Assume that $\pi$ is inner faithful.
Let $I \subset {\rm Ker}(\pi)$ be a bi-ideal. Then $J = \sum_{k \in \mathbb N} S^k(I)$
is a coideal stable under the antipode with $J \subset {\rm Ker}(\pi)$. Let $J'$ be the ideal
generated by $J$: this is a bi-ideal with $S(J')\subset J'$, and hence  $J'$ is a Hopf ideal
contained in ${\rm Ker}(\pi)$. Hence by the previous proposition we have
$J=J'=(0)$. The converse assertion is immediate.
\end{proof}

We also have a characterization of the Hopf image using inner faithfulness,
which can be useful in some contexts (for example see the proof of Proposition 3.3).

\begin{proposition}
  Let $\pi : H \longrightarrow A$ be a representation
of a Hopf algebra $H$ on an algebra $A$, and let 
$(H_\pi, p, \tilde{\pi})$ be the Hopf image of $\pi$. 
Then $\tilde{\pi}$ is inner faithful. Conversely if $(L,q, \varphi)$ is a factorization of $\pi$
such that $\varphi$ is inner faithful, then $L \simeq H_\pi$.
\end{proposition}

\begin{proof}
Let $(L,q,\varphi)$ be a factorization of $\tilde{\pi} : H_\pi \longrightarrow A$: we have 
$\varphi \circ q = \tilde{\pi}$, and hence $\varphi \circ (q \circ p) = \pi$.
Thus there exists a Hopf algebra map $f : L \longrightarrow H_\pi$
such that $f \circ q \circ p = p$. Hence $f \circ q = {\rm id}_{H_\pi}$, 
$q$ is injective and is an isomorphism. By Proposition 2.8 
we conclude that $\tilde{\pi}$ is inner faithful.
 
Now let
$(L,q, \varphi)$ be a factorization of $\pi$ with $\varphi$ inner faithful.
By the universal property of the Hopf image, there exists a surjective Hopf algebra map
$f : L \longrightarrow H_\pi$ such that $\tilde{\pi} \circ f = \varphi$.
But then $(H_\pi, f, \tilde{\pi})$ is a factorization of $\varphi$, and Proposition 2.8 ensures
that $f$ is an isomorphism. 
\end{proof}

When considering Hopf algebra maps, Hopf images are usual images and inner faithfulness is equivalent
to faithfulness.

\begin{proposition}
 Let $H$ and $A$ be Hopf algebras and $\pi : H \longrightarrow A$ be a Hopf algebra map.
 Then the Hopf image of $\pi$ is $\pi(H)$, and
$\pi $ is inner faithful if and only if it is faithful. 
\end{proposition}

 The first assertion is immediate, as well as the second one, since
 ${\rm Ker}(\pi)$ is a Hopf ideal.
The following result is also immediate using the universal property of the Hopf image. 

\begin{proposition}
 Let $\pi : H \longrightarrow A$ be a representation
of a Hopf algebra $H$ on an algebra $A$ and let 
$\theta : A \longrightarrow B$ be an algebra isomorphism. Then we have a Hopf algebra isomorphism
$H_{\theta \circ \pi} \simeq H_\pi$. In particular $\theta \circ \pi$ is inner faithful
if and only if $\pi$ is inner faithful.
\end{proposition}

Of course one cannot expect that an inner faithful representation $H \longrightarrow A$
will transmit all the algebra properties
of the  algebra $A$ to the algebra $H$. This is true however for commutativity.

\begin{proposition}
 Let $\pi : H \longrightarrow A$ be a representation
of a Hopf algebra $H$ on a commutative algebra $A$. Then the Hopf image $H_\pi$
is a commutative algebra, and hence if $\pi$ is inner faithful, then $H$ is also commutative.
\end{proposition}

\begin{proof}
 Since $A$ is commutative, the ideal of $H$ generated by the commutators of elements of $H$
 is contained in ${\rm Ker}(\pi)$. But the commutator ideal  is also a Hopf ideal, so 
 is contained in $I_\pi$, and hence $H_\pi$ is commutative.
\end{proof}

\subsection{Hopf $*$-algebras.}
In this subsection we assume that $k=\C$, and we consider
Hopf $*$-algebras. In this case the construction of the previous subsection has to be slightly modified.

Let us begin by recalling the appropriate language.
First, a Hopf $*$-algebra is a Hopf algebra $H$ which is a $*$-algebra and 
such that the comultiplication $\Delta : H \longrightarrow H \otimes H$ is a $*$-algebra map. It follows that the counit is a $*$-morphism and that
 the antipode is bijective, and its inverse
 satisfies $S^{-1}(x) =  S (x^*)^*$. 
 A $*$-representation of $H$ on a $*$-algebra is a 
 $*$-algebra map $H \longrightarrow A$.

 The formulation of the problem of the existence of a Hopf $*$-image is the same as in the introduction, adding ``$*$'' where needed. In this framework, Theorem 2.1 has the following form.

 \begin{theorem}
Let $\pi : H \longrightarrow A$ be a $*$-representation
of a Hopf $*$-algebra $H$ on an $*$-algebra $A$. Then $\pi$
has a Hopf $*$-image:
there exists a Hopf $*$-algebra $H_\pi$ together
with a surjective Hopf $*$-algebra map $p : H \longrightarrow H_\pi$ and 
 a $*$-representation $\tilde{\pi} : H_\pi \longrightarrow A$, such that
$\pi = \tilde{\pi} \circ p$, and such that if $(L,q, \varphi)$ is another 
$*$-factorization of $\pi$, there exists a unique
Hopf $*$-algebra map $f : L \longrightarrow H_\pi$ such that
$f \circ q = p$ and $\tilde{\pi} \circ f = \varphi$.  
$$\xymatrix{
H \ar[rr]^{\pi} \ar[dr]_{p} \ar@/_/[ddr]_q & & A  \\
& H_\pi \ar[ur]_{\tilde{\pi}} & \\
& L \ar@{-->}[u]_f \ar@/_/[uur]_{\varphi}& 
}$$
\end{theorem}

\begin{proof} Similarly to the proof of Theorem 2.1, one constructs
$I_\pi^+$, 
the largest Hopf $*$-ideal contained in ${\rm Ker}(\pi)$.
The explicit construction of $I_\pi^+$ is as follows.
Let $F^+$ be the free monoid generated by the set $\mathbb Z$. Similarly
to the construction in the previous subsection, to any $g \in F^+$, we associate an algebra 
$A^g$ and a representation $\pi^g : H \longrightarrow A^g$, using the inverse 
of the antipode when needed. Then, similarly to the proof of Proposition 2.2, we show that
$$I_\pi^+ = \bigcap_{g \in F^+} {\rm Ker}(\pi^g)\subset H$$
 The details are left to the reader (to show that $\bigcap_{g \in F^+} {\rm Ker}(\pi^g)$ is $*$-stable, one uses the formula $S\circ * = * \circ S^{-1})$. 
\end{proof}

\begin{remark}
 {\rm 
 We have used the same notation $H_\pi$ for the Hopf $*$-image of a $*$-representation
 $\pi : H \longrightarrow A$ and its Hopf image. In general it seems possible that
 two notions might not coincide, although we do not have an explicit example. 
 This should cause no confusion: unless
 specifically notified, when we have such a $*$-representation, the notation $H_\pi$ will always denote the Hopf $*$-image.
   
  We will say that a $*$-representation $\pi : H \longrightarrow A$
of a Hopf $*$-algebra $H$ on an $*$-algebra $A$ is inner faithful if 
${\rm Ker}(\pi)$ does not contain any non zero Hopf $*$-ideal
(similarly to Proposition 2.8
there are several equivalent characterizations). Again it seems to be possible  
that such a $*$-representation  be inner faithful as a $*$-representation,
but not as a representation. However the two notions coincide when $S^2= {\rm id}_H$
(or more generally if some power of $S^2$ is an inner automorphism).
   }
\end{remark}

\begin{remark}
 {\rm 
 The proof of Theorem 2.14 also shows that Hopf images exist in the category of Hopf algebras
 having a bijective antipode (more precisely, the Hopf ideal $I^+_\pi$ is the largest Hopf ideal 
 with $S(I_\pi^+) = I_\pi^+$ contained in ${\rm Ker}(\pi)$), and 
 the Hopf image in this category coincides with the Hopf $*$-image.
  Note that the construction of Theorem 2.1 
 does not give the Hopf image in the category of Hopf algebras
 with bijective antipode, since there exist quotients of 
 Hopf algebras with bijective antipode that do not have a bijective antipode: see e.g.
 \cite{sc00}.
   }
\end{remark}

We now turn to \textbf{compact Hopf algebras}: these are the Hopf $*$-algebras 
having all their finite-dimensional comodules equivalent 
to unitary ones (see \cite{ks}, these are called CQG algebras there), and
are the algebras of representative functions on compact quantum groups.

\begin{theorem}
 Let $\pi : H \longrightarrow A$ be a $*$-representation
of a compact Hopf algebra $H$ on a $*$-algebra $A$.
Then the Hopf $*$-image $H_\pi$ is a compact Hopf algebra, and hence 
Hopf $*$-images exist in the category of compact Hopf algebras.
\end{theorem}

\begin{proof} 
It is clear from Theorem 27, Section 11.3 in \cite{ks}, that
a quotient of a compact Hopf algebra is again a compact Hopf algebra.
Thus since we have a surjective Hopf $*$-algebra map $H \longrightarrow H_\pi$, we conclude
that $H_\pi$ is a compact Hopf algebra.  
\end{proof}

\section{Classical examples: group algebras and Lie algebras}

In this section we get back to the motivating examples for the notion of inner faithfulness:
group algebras and universal enveloping algebras of Lie algebras.

\subsection{Group algebras} The following result
is announced in the introduction. Its origin goes back to Proposition 2.2 in \cite{ba98}.

\begin{proposition}
 Let $ \Gamma$ be group, let $A$ be an algebra and let $\pi : k[\Gamma] \longrightarrow A$
 be an algebra map. Then we have 
 $$k[\Gamma]_\pi \simeq k[\pi(\Gamma)] \simeq k[\Gamma/{\rm Ker}(\pi_{|\Gamma})]$$
 In particular $\pi$ is inner faithful if and only if
 the group morphism $\pi_{|\Gamma} : \Gamma \longrightarrow A^\times$ is injective (faithful).
\end{proposition}

\begin{proof}
 Let $(L,q, \varphi)$ be a factorization of $\pi$. Then 
 $L$ is the group algebra $k[q(\Gamma)]$, and we have a group factorization
 $$\xymatrix{
\Gamma \ar[rr]^{\pi_{|\Gamma}} \ar[dr]_{q_{|\Gamma}} & & A^\times \\
& q(\Gamma) \ar[ur]_{\varphi_{|q(\Gamma)}} &
}$$
Hence we have $q(\Gamma) \subset \pi(\Gamma)$, and this induces the appropriate
Hopf algebra map $k[q(\Gamma)] =L \longrightarrow k[\pi(\Gamma)]$.  
\end{proof}

The above result motivates the following definition.

\begin{definition}
 A Hopf algebra $H$ is said to be \textbf{inner linear}
 if there exists an inner faithful representation  $\pi : H \longrightarrow A$
 into a finite-dimensional algebra $A$.   
\end{definition}

Indeed the group algebra $k[\Gamma]$ is inner linear if and only the 
group $\Gamma$ is linear (over $k$). It is clear that a Hopf algebra 
is inner linear if and only if it contains an ideal of finite codimension
that does not contain non-zero Hopf ideals.

\subsection{Lie algebras}
Let $\mathfrak g$ be a Lie algebra, and let $U(\mathfrak g)$ 
be its universal enveloping algebra. Algebra maps $U(\mathfrak g) \longrightarrow A$ correspond
exactly to Lie algebra maps $\mathfrak g \longrightarrow \mathfrak{gl}(A)$, where 
$\mathfrak{gl}(A)$ is the Lie algebra associated to $A$: as a vector space 
$\mathfrak{gl}(A)=A$ and the Lie bracket is defined by $[a,b] = ab-ba$, $\forall a,b \in A$. 
The Hopf image is described as follows in characteristic zero.

\begin{proposition}
 Let $ \mathfrak g$ be a Lie algebra over a characteristic zero field, let $A$ be an algebra and let $\pi : U(\mathfrak g) \longrightarrow A$
 be an algebra map. Then we have 
 $$U(\mathfrak g)_\pi \simeq U(\pi(\mathfrak g)) \simeq U(\mathfrak g/{\rm Ker}(\pi_{|\mathfrak g}))$$
 In particular $\pi$ is inner faithful if and only if
 the Lie algebra map $\pi_{|\mathfrak g} : \mathfrak g \longrightarrow \mathfrak{gl}(A)$ is injective (faithful). Thus the Hopf algebra $U(\mathfrak g)$ is inner linear if and only if
 $\mathfrak g$ is finite-dimensional.
\end{proposition}

\begin{proof}
 Let us first show if $\pi_{|\mathfrak g} : \mathfrak g \longrightarrow \mathfrak{gl}(A)$ is injective,
 then $\pi : U(\mathfrak g) \longrightarrow A$ is inner faithful. 
 The space of primitive elements $\mathcal P(U(\mathfrak g))$ equals $\mathfrak g$
 by the characteristic zero assumption and $1$ is the only group-like in $U(\mathfrak g)$
 (see e.g. Proposition 5.5.3 in \cite{mo}). Thus the canonical map $p : U(\mathfrak g) \longrightarrow
 U(\mathfrak g)_\pi$ is injective on primitive elements, and since the Hopf algebra $U(\mathfrak g)$
 is pointed, we use Corollary 5.4.7 in \cite{mo} to conclude that $p$ is injective, and
 hence $\pi$ is inner faithful.
 
 In general we have a Hopf algebra factorization
 $$\xymatrix{
U(\mathfrak g) \ar[rr]^{\pi} \ar[dr] & & A \\
& U(\pi(\mathfrak g)) \ar[ur] &
}$$
The algebra map on the right is inner faithful by the previous discussion, and hence 
by Proposition 2.9 we have the 
announced isomorphism
$U(\mathfrak g)_\pi \simeq U(\pi(\mathfrak g)) \simeq U(\mathfrak g/{\rm Ker}(\pi_{|\mathfrak g}))$.
Finally if $\pi$ is inner faithful, it induces an isomorphism $U(\mathfrak g) \simeq U(\pi(\mathfrak g))$, which
induces an isomorphism   between the Lie algebras of primitive elements, and hence
between $\mathfrak g$ and $\pi(\mathfrak g)$: $\pi_{|\mathfrak g}$ is injective.

If $U(\mathfrak g)$ is inner linear, then $\mathfrak g$ is finite-dimensional by the previous discussion. The converse follows from Ado's Theorem: a finite dimensional
Lie algebra $\mathfrak g$ has a faithful finite-dimensional representation
$\mathfrak g \longrightarrow \mathfrak{gl}_n(k)$ for some $n$. 
\end{proof}

In fact the main argument we have used for the proof this proposition, Corollary 
5.4.7 in \cite{mo}, a result due independently  to Takeuchi and Radford,
 is valid for arbitrary pointed Hopf algebras.
 We use it in a similar manner in the next section to get an inner faithfulness
 criterion for representations  of arbitrary pointed Hopf algebras.

\section{Pointed Hopf algebras}

After the motivating examples of group algebras and enveloping
algebras of Lie algebras, the next natural step is the study
of pointed Hopf algebras. These have been much studied in recent years
(see e.g. \cite{ks} for quantized enveloping algebras of Lie algebras and \cite{as}
for the finite-dimensional case).
In this section we study the inner faithfulness of a representation
of a pointed Hopf algebra. The criterion that we give unifies those 
given in the previous section.

Recall that a Hopf algebra $H$ is said to be pointed if all its simple 
comodules are one-dimensional, hence corresponding to group-like elements.
The group of group-like elements of $H$ is denoted by Gr$(H)$.
An element $x \in H$ is said to be $(g,h)$-primitive for 
$g,h \in {\rm Gr}(H)$ if
$$\Delta(x) = g \otimes x + x \otimes h$$
The space of $(g,h)$-primitive elements is denoted $\mathcal P_{g,h}(H)$.

\begin{theorem}
Let $\pi : H \longrightarrow A$ be a representation
of a pointed Hopf algebra $H$ on an algebra $A$.
The following assertions are equivalent.
\begin{enumerate}
 \item $\pi$ is inner faithful.
 \item $\forall g\in {\rm Gr}(H)$, the 
restriction $\pi_{|\mathcal P_{g,1}(H)} : \mathcal P_{g,1}(H) \longrightarrow A$ 
is injective.
\item $\forall g \in {\rm Gr}(H)$, the 
restriction $\pi_{|\mathcal P_{1,g}(H)} : \mathcal P_{1,g}(H) \longrightarrow A$ 
is injective.
\end{enumerate}
\end{theorem}

\begin{proof}
We begin by showing that $\pi$ is inner faithful if and only if
for all $g,h\in {\rm Gr}(H)$, the 
restriction $\pi_{|\mathcal P_{g,h}(H)} : \mathcal P_{g,h}(H) \longrightarrow A$ 
is injective.
Assume that $\pi$ is inner faithful. Let $x \in  \mathcal P_{g,h}(H)$
be such that $\pi(x) =0$ and let $I$ be the ideal of $H$
generated by $x$. It is clear from the identities 
$$\Delta (x)= g \otimes x + x \otimes h \quad {\rm and} \quad S(x) =-g^{-1}xh^{-1}$$
that $I$ is a Hopf ideal, and since $I$ is contained in ${\rm Ker}(\pi)$,
we get $I=0$ and $x=0$.  

Conversely, assume that $\pi$ is injective on each space of primitives. Then
so is $p : H \longrightarrow H_\pi$, and hence by Corollary 5.4.7
in \cite{mo}, we conclude that $p$ is an isomorphism and that
$\pi$ is inner faithful.

For $x \in \mathcal P_{g,h}(H)$, we have $g^{-1}x \in \mathcal P_{1,g^{-1}h}(H)$
and $h^{-1}x  \in \mathcal P_{h^{-1}g,1}(H)$. Hence if (2) or (3) holds, then
$\pi$ is injective when restricted to each space of primitives, and the previous discussion
shows that $\pi$ is inner faithful.
\end{proof}

The above criterion applies in particular to the quantized universal enveloping algebra
$U_q(\mathfrak g)$ of a semisimple Lie algebra $\mathfrak g$, for which 
descriptions of skew-primitives are available \cite{cm,mu}.
As a simple illustration, let us have a look at 
the vector representation of $U_q(\mathfrak{sl}_2(\C))$.

\begin{theorem}
 Let $\pi : U_q(\mathfrak{sl}_2(\C)) \longrightarrow M_2(\C)$ be 
 the 2-dimensional vector representation of $U_q(\mathfrak{sl}_2(\C))$.
 Then $\pi$ is inner faithful if and only if
 $q$ is not a root of unity, and hence $U_q(\mathfrak{sl}_2(\C))$
 is inner linear if $q$ is not a root of unity.
\end{theorem}

\begin{proof}
 Recall first that for $q \in \C$, $q^2\not = 1$ and $q \not = 0$, 
 the algebra $U_q=U_q(\mathfrak{sl}_2(\C))$ is 
 presented by the generators $E$, $F$, $K$, $K^{-1}$, submitted to the relations
 $$KK^{-1} = 1 = K^{-1}K, \ KE= q^2EK, \ KF = q^{-2}FK, \ 
 [E,F] = \frac{K-K^{-1}}{q-q^{-1}}$$
 and its comultiplication is defined by
 $$\Delta(E) = 1 \otimes E + E \otimes K, \quad  \Delta(F) = K^{-1} \otimes F + F \otimes 1, \quad
 \Delta(K) = K \otimes K$$
 The vector representation $\pi : U_q \longrightarrow M_2(\C)$
 is defined by
 $$\pi(K) = \begin{pmatrix} q & 0 \\ 0 & q^{-1}\end{pmatrix}, \ 
 \pi(E) = \begin{pmatrix} 0 & 1 \\ 0 & 0\end{pmatrix}, \ 
 \pi(K) = \begin{pmatrix} 0 & 0 \\ 1 & 0\end{pmatrix} $$
 We have Gr$(U_q) =\{K^n, \ n \in \mathbb Z\}\simeq \mathbb Z$.
 If $q$ is not a root of unity we have
 $$\mathcal P_{1,K}(U_q) = \C(1-K) \oplus \C E \oplus \C FK, \  {\rm and} \ 
 {\rm for} \  m \not = 1, \ \mathcal P_{1,K^m}(U_q) = \C(1-K^m)$$
 It is a direct computation to check that the conditions of Theorem 4.1 are fulfilled,
 and hence $\pi$ is inner faithful.
 If $q$ is a root of unity, then $\pi$ is not injective 
 on the group-like elements and hence is not inner faithful. 
\end{proof}

\begin{remark}
 {\rm If $q$ is a root of unity, one can show that the Hopf image
 of the above representation $\pi$ is the (finite-dimensional) Hopf algebra
 $u_q(\mathfrak sl_2(\C))$.}
\end{remark}

\section{Function algebras}

In this section we study Hopf images for various function algebras: polynomial
functions on algebraic groups and representative functions on  compact groups.
The idea for the computation of the Hopf image goes back to \cite{ba98}, but 
we are a little bit more general here.  
These simple examples are already interesting for testing
the possibility of generalizing representation theoretic
properties of discrete group algebras to arbitrary Hopf algebras.

\begin{proposition}
 Let $G$ be a linear algebraic group over an algebraically closed field of characteristic zero
 and let $\pi : \mathcal O(G) \longrightarrow A$ be a representation
 on an algebra $A$. Assume that the algebra $\pi(A)$ is reduced, so that
 $\pi(A) \simeq \mathcal O(X)$ for an affine algebraic set $X$
 and that the algebra map $\mathcal O(G) \longrightarrow \pi(A) \simeq \mathcal O(X)$
 is induced by a polynomial map $\nu: X \longrightarrow G$.
 Then $\mathcal O(G)_\pi \simeq \mathcal O(\overline{\langle \nu(X)\rangle})$, where
  $\overline{\langle \nu(X)\rangle}$ is the Zariski closure in $G$ of the subgroup generated
  by $\nu(X)$. 
\end{proposition}

\begin{proof}
By the assumption and Proposition 2.12 we may assume that 
$\pi : \mathcal O(G) \longrightarrow \mathcal O(X)$ is induced by an injective
polynomial map $\nu : X \longrightarrow G$.
The injections $X \rightarrow  \overline{\langle \nu(X)\rangle} \subset G$
yield  a factorization of $\pi$.
Now let $(L, q, \varphi)$ be a factorization of $\pi$. 
Then $L$ is finitely generated and is reduced by Cartier's theorem
  ($k$ has characteristic zero),  hence we can assume
  that $L= \mathcal O(H)$ for a linear algebraic group $H$, and that 
  $q$ and $\rho$ are induced by polynomial maps $X \rightarrow H \subset G$
  whose composition is $\nu$. Hence $\overline{\langle \nu(X)\rangle} \subset H$, and this gives
  the required Hopf algebra map 
  $\mathcal O(H) \longrightarrow \mathcal O(\overline{\langle \nu(X)\rangle})$.
\end{proof}

\begin{example}{\rm
Let $G$ be a linear algebraic group over an algebraically closed field of characteristic zero
 and let $g_1, \ldots ,g_n \in G$.
 The algebra map
 \begin{align*}
  \mathcal O(G) &\longrightarrow k^n \\
  f &\longmapsto (f(g_1), \ldots , f(g_n))
 \end{align*}
has $\mathcal O(\overline{\langle g_1, \ldots ,g_n\rangle})$ as Hopf image
and is inner faithful if and only if $G = \overline{\langle g_1, \ldots ,g_n\rangle}$.
} 
\end{example}

The compact group case is worked out similarly.
Let $G$ be a compact group.
The Hopf $*$-algebra of representative functions on $G$ is denoted by $\mathcal R(G)$:
this a dense $*$-subalgebra of C$(G)$ by the Peter-Weyl theorem, and moreover C$(G)$ is the enveloping 
$\c^*$-algebra of $\mathcal R(G)$.

\begin{proposition}
 Let $G$ be a compact group and let
 $\pi : \mathcal R(G) \longrightarrow A$
 be a $*$-representation on a ${\rm C^*}$-algebra $A$.
 Extend $\pi$ to a $*$-algebra map $\pi^+ : C(G) \longrightarrow A$, and 
 let $X$ be the spectrum of $\pi^+(C(G))$, so that
 $\pi^+(C(G)) \simeq C(X)$ and that the induced 
 ${\rm C^*}$-algebra map $C(G) \longrightarrow \pi^+(C(G))\simeq C(X)$ comes from a 
 continuous map $\nu : X \longrightarrow G$. 
 Then $\mathcal R(G)_\pi \simeq \mathcal R(\overline{\langle \nu(X)\rangle})$, where
  $\overline{\langle \nu(X)\rangle}$ is the closure in $G$ of the subgroup generated
  by $\nu(X)$. 
\end{proposition}

\begin{proof}
The proof is essentially the same as the one of the previous proposition,
because a quotient of a compact Hopf algebra is itself compact, and
by the Hopf algebra version of the Tannaka-Krein theorem, a commutative 
compact Hopf algebra is the algebra of representative functions on a unique compact group. 
\end{proof}

\begin{theorem}
 Let $G$ be a compact Lie group. Then there exists an inner faithful 
 $*$-representation $\pi : \mathcal R(G) \longrightarrow \C^n$ for some $n \in \mathbb N^*$.
\end{theorem}

\begin{proof}
If $G$ is connected, it is known \cite{au} that there exists $x, y \in G$ such
$G = \overline{\langle x , y\rangle}$. Similarly to Example 5.2, this gives a 
$*$-representation $\mathcal R(G) \longrightarrow \C^2$ which is inner faithful by the previous proposition. If $G$ is not connected, let $G_0$ be the neutral component of $G$: 
the group $G/G_0$ is finite and hence $G$ has a family of $2[G:G_0]$  topological
generators, which again gives an inner faithful representation on $\C^n$, with
$n = 2 [G:G_0]$.
\end{proof}

\begin{corollary}
 Let $G$ be a complex linear algebraic group. Then there exists an inner faithful 
 representation $\pi : \mathcal O(G) \longrightarrow \C^n$ for some $n \in \mathbb N^*$.
\end{corollary}

\begin{proof}
If $G$ is reductive,
 this follows from the previous result since the Hopf algebra $\mathcal O(G)$ is the 
 algebra of representative functions on a maximal compact subgroup of $G$. 
 Then the Levi decomposition of $G$ as a semi-direct product of a reductive
 group and a unipotent one reduces the problem to the case when $G$ is unipotent.
 If $G$ is additive group of $\C$, the result is well-known. 
 Then the case of the group of unipotent matrices is proved by induction,  and
 the general case follows by using the representation of a unipotent group as a group of unipotent
 matrices.  
\end{proof}

We end the section by showing
that the simple example of function algebras
shows that it is not possible to deduce 
inner faithfulness of a representation of a cosemisimple Hopf algebra
by only studying its restriction to characters.

Let us assume that $k$ is algebraically closed.
Recall that a cosemisimple Hopf algebra is a Hopf algebra $H$
whose category of comodules is semisimple. This is equivalent to say that $H$ has a 
Peter-Weyl decomposition
$$H = \bigoplus_{\lambda \in \Lambda} H_\lambda$$
where $\Lambda$ is the set of simple $H$-comodules and for $\lambda \in \Lambda$, $H_\lambda$
is the comatrix coalgebra of corresponding coefficients. Let 
$d_\lambda$ be the dimension of the simple $H$-comodule corresponding to $\lambda$, and let $\chi_\lambda \in H_\lambda$ be the corresponding character.

The group algebra case corresponds to when $d_\lambda=1$ for any $\lambda$, and
then the characters correspond to the group-like elements.
Since for group algebras inner faithfulness can be detected
by only studying the restriction of a representation to group-like elements, it would be natural
to hope that in the general case, the restriction to characters would lead to the same
information. The following simple example shows that this is not true.   

\begin{example}{\rm 
Let us consider $S_4$, the symmetric group on $4$ letters, and let
$\C(S_4)=\C[S_4]^*$ be the cosemisimple Hopf algebra of (complex) functions on $S_4$.
The group $S_4$ is generated by the elements $\tau_1=(1,2)$, $\tau_2=(2,3)$ and $\tau_3=(3,4)$, and hence 
we have an inner faithful representation
\begin{align*}
  \C(S_4) &\longrightarrow \C^3 \\
  f &\longmapsto (f(\tau_1), f(\tau_2), f(\tau_3))
 \end{align*}
The elements   $\tau_1$, $\tau_2$ and $\tau_3$  are all conjugate in $S_4$, so 
this representation is not injective on characters. Thus an inner faithful representation
is not necessarily injective on characters.

Conversely the representation
\begin{align*}
  \C(S_4) &\longrightarrow \C^2 \\
  f &\longmapsto (f(\tau_1), f((1,2,3)))
 \end{align*}
is injective on characters (as one easily checks on the character table of $S_4$)
but is not inner faithful since the elements $\tau_1$ and $(1,2,3)$ do not generate
$S_4$.

As a conclusion, it seems that there is no link between the inner faithfulness
of a representation and its injectivity on characters.
}
\end{example}

\section{Twistings}

In this  section we study the behaviour of Hopf images under
various deformation procedures such as Drinfeld's twisting,
or the dual operation, often called 2-cocycle deformation, that
we will call here cotwisting for simplicity.   

\subsection{Twisting} We begin with the
twisting operation, in the sense of Drinfeld \cite{dr}.
More exactly the definitions we use are taken or adapted from 
\cite{ev,ni}.

Let $H$ be a Hopf algebra
and let $\Omega$ be an invertible element in $H\otimes H$.
Consider the linear maps
 $\delta_\Omega, \Delta_\Omega: H \longrightarrow H \otimes H$
  defined respectively by
 $$\delta_\Omega(x) = \Omega \Delta(x), \quad \Delta_\Omega(x) = \Omega \Delta(x) \Omega^{-1}$$
 We say that $\Omega$ is a \textbf{twist} on $H$ if $(H, \delta_\Omega,\varepsilon)$
is a coalgebra. The element $u = m \circ ({\rm id_H} \otimes S)(\Omega)$
is then invertible in $H$, and $H_\Omega=(H,m,u,\Delta_\Omega, \varepsilon, S_u)$
is a Hopf algebra, where $S_u : H\longrightarrow H$ is defined by $S_u(x)= uS(x)u^{-1}$.
 
 We say that $\Omega$ is a \textbf{pseudo-twist} on $H$  if $(\varepsilon\otimes {\rm id})(\Omega) = 1 = 
({\rm id} \otimes \varepsilon)(\Omega)$, if $(H, \Delta_\Omega,\varepsilon)$
is a coalgebra, and if there exists an invertible element
$u \in H$ such that $S_u : H \longrightarrow H$, defined by $S_u(x) =u S(x)u^{-1}$,  is an 
antipode for the bialgebra $(H,m,u,\Delta_\Omega, \varepsilon)$, so that
$H_\Omega =  (H,m,u,\Delta_\Omega, \varepsilon, S_u)$ is a Hopf algebra.

A Hopf algebra having the form $H_\Omega$ for some twist (resp. pseudo-twist) $\Omega$
on $H$ is said to be a \textbf{twist} (resp. \textbf{pseudo-twist}) of $H$.

The following lemma gives some basic properties of twisting, probably well-known.
The direct verification is left to the reader.

\begin{lemma}
Let $\Omega$ be a pseudo-twist on a Hopf algebra $H$.
\begin{enumerate}
\item Let  $f : H \longrightarrow L$ be a surjective Hopf algebra map. Then 
$f^*(\Omega) = (f \otimes f)(\Omega)$ is also a pseudo-twist for $L$.
\item  The
Hopf ideals in $H$ are exactly the Hopf ideals in $H_\Omega$. 
\end{enumerate}
\end{lemma}

\begin{theorem}
 Let $H$ be a Hopf algebra, let $\Omega$ be a pseudo-twist on $H$ and let $\pi : H \longrightarrow A$ be a representation on an algebra $A$, that we also view as a representation
 $\pi : H_\Omega \longrightarrow A$. We have $(H_\Omega)_\pi =
 (H_\pi)_{p^*(\Omega)}$, where $p : H \longrightarrow H_\pi$ is the canonical projection,
 and $\pi$ is inner faithful as a representation of $H$ if and only if it is inner faithful as a 
 representation of $H_\Omega$.
\end{theorem}

\begin{proof}
 Since by the previous lemma the Hopf ideals of $H$ and of $H_\Omega$ are the same, this
 is also true for the Hopf ideals contained in ${\rm Ker}(\pi)$, and the largest one is the 
 same. So the defining Hopf ideal of $(H_\Omega)_\pi$ is the Hopf ideal $I_\pi$ constructed in Section 2, so that
 $(H_\Omega)_\pi = H_\Omega/ I_\pi$, and since $(H/I_\pi)_{p^*}(\Omega) = H_\Omega/I_\pi$, we have the claimed result. The last assertion is also immediate after the previous discussion.
\end{proof}

It follows that the twist of an inner linear Hopf algebra 
is still inner linear.
The reader will find several examples of twisted Hopf algebras in \cite{ev,ni}, for example, for which
the above theorem furnishes inner faithful representations of the new Hopf algebra from the old one. Some of the Hopf algebras considered in Section 9 are also obtained by twisting.

\subsection{Cotwisting} 
We now deal with the dual operation to twisting, that
we call cotwisting. 
The situation is more complicated here, because we deform the product rather
than the coproduct, and so a representation of the given Hopf algebra
is not a representation of the deformed one.

Let us recall the basic vocabulary, which is dual to the one of the previous subsection.
We only consider the cotwist case.

Let $H = H$ be a Hopf algebra.
We use Sweedler's notation
$\Delta(x) = x_{1} \otimes x_{2}$. Recall (see e.g. \cite{do})
that a \textbf{cotwist} (=2-cocycle) is a convolution invertible linear map
$\sigma : H \otimes H \longrightarrow \c$ satisfying
$$\sigma(x_{1}, x_{2}) \sigma(x_{2}y_{2},z) =
\sigma(y_{1},z_{1}) \sigma(x,y_{2} z_{2})$$
and $\sigma(x,1) = \sigma(1,x) = \varepsilon(x)$, for $x,y,z \in H$.

Following \cite{do} and \cite{sc1}, we associate various algebras to
a cotwist. First consider the algebra 
$_{\sigma} \! H$. As a vector space we have $_{\sigma} \! H = H$ and the product
of $_{\sigma}H$ is defined to be
$$\{x\}  \{y\} = \sigma(x_{1}, y_{1}) \{x_{2} y_{2}\}, 
\quad x,y \in H,$$
where an element $x \in H$ is denoted $\{x\}$, when viewed as an element 
of $_{\sigma} \!H$.

We also have the algebra $H_{\sigma^{-1}}$, where $\sigma^{-1}$ denotes the convolution inverse of $\sigma$. As a vector space we have
$H_{\sigma^{-1}} = H$ and the product of 
$H_{\sigma^{-1}}$ is defined to be
$$\langle x \rangle \langle y \rangle = \sigma^{-1}(x_{2}, y_{2}) \langle x_{1} y_{1} \rangle, 
\quad x,y \in H,$$
where an element $x \in H$ is denoted $\langle x \rangle$, when viewed as an element 
of $H_{\sigma^{-1}}$. The cocycle condition ensures that $_{\sigma} \! H$
and $H_{\sigma^{-1}}$ are associative algebras with $1$ as a unit.

Finally we have the 
Hopf algebra $H^{\sigma} = {_{\sigma} \! H}\!_{\sigma^{-1}}$.  
As a coalgebra $H^{\sigma} 
= H$. The product of $H^{\sigma}$ is defined to be
$$[x] [y]= \sigma(x_{1}, y_{1})
\sigma^{-1}(x_{3}, y_{3}) [x_{2} y_{2}], 
\quad x,y \in H,$$
where an element $x \in H$ is denoted $[x]$, when viewed as an element 
of $H^{\sigma}$, 
and we have the following formula for the antipode of 
$H^{\sigma}$:
$$S^{\sigma}([x]) = \sigma(x_{1}, S(x_{2}))
\sigma^{-1}(S(x_{4}), x_{5})[ S(x_{3})].$$
A Hopf algebra having the form $H^\sigma$ for some cotwist  $\sigma$
on $H$ is said to be a \textbf{cotwist} of $H$.

Very often cotwists are induced by simpler quotient Hopf algebras.
More precisely let $\pi : H \to K$ be a Hopf algebra surjection and let
$\sigma : K \otimes K \to \C$ be a cotwist on $K$.
Then $\sigma_{\pi} = \sigma \circ (\pi \otimes \pi) : H \otimes H \to \C$ is a cotwist.

We prove the following result.

\begin{theo}\label{cotwist}
 Let $H$ be a Hopf algebra  and let $\sigma : H \otimes H \to k$
be a cotwist induced by a finite-dimensional quotient Hopf algebra
of $H$. Assume that $S^2$ is an inner automorphism of $H$. If $H$ is inner linear, then
$H^\sigma$ is also inner linear.
\end{theo}

The proof will be a consequence of the following technical result.

\begin{theorem}
 Let $\theta : H \longrightarrow A \otimes L$ be an  algebra map, where $H,L$ are Hopf algebras and $A$ is an algebra, and let $\varphi : L \longrightarrow B$
 be an inner faithful representation of $L$ on an algebra $B$. 
 Assume that there exists $\psi \in A^*$ such
 that $(\psi \otimes {\rm id}_L) \circ \theta : H \longrightarrow L$ is an injective coalgebra map.
 Assume moreover that one of the following conditions holds.
 \begin{enumerate}
  \item $S_L \circ ((\psi \otimes {\rm id}_L)\circ \theta) = ((\psi \otimes {\rm id}_L) \circ \theta)\circ S_H$.
  \item $S_L({\rm Ker}(\varphi)) \subset {\rm Ker}(\varphi)$.
 \end{enumerate}
Then the representation $({\rm id}_A \otimes \varphi) \circ \theta : H \longrightarrow A \otimes B$
is inner faithful.
\end{theorem}

\begin{proof}
 Let $I \subset {\rm Ker}(({\rm id}_A \otimes \varphi) \circ \theta)$ be a Hopf ideal. Let
 $J = (\psi \otimes {\rm id}_L)(\theta(I))\subset L$: this is a coideal since 
 $(\psi \otimes {\rm id}_L) \circ \theta$ is a coalgebra map. We have 
 $$\varphi(J) = (\psi \otimes \varphi)(\theta(I)) = (\psi \otimes {\rm id}_L) ( ({\rm id}_A \otimes \varphi)(\theta(I))=0$$
 and hence $J \subset {\rm Ker}(\varphi)$. 
 Let $J^+$ be the ideal of $L$ generated by $J$: we have $J^+ \subset {\rm Ker}(\varphi)$ and
 $J^+$ is a bi-ideal in $L$.
 
 Assume that condition (1) holds. Then $J$ is $S_L$-stable and hence 
 so is $J^+$, which is a Hopf ideal. Since $\varphi$ is inner faithful, we have
 $J^+=J= (0)$, and we conclude that $I=(0)$ by the injectivity of $(\psi \otimes {\rm id}_L) \circ \theta$.
 
 Assume now that condition (2) holds. Then $J^+=(0)=J$ by Proposition 2.8, and $I=(0)$, again by
 the injectivity of $(\psi \otimes {\rm id}_L) \circ \theta$.
\end{proof}

\begin{proof}[Proof of Theorem 6.3]
By the assumption, there is a Hopf algebra surjection
 $\pi : H \to K$  onto  a finite dimensional Hopf algebra $K$ and a cotwist
$\tau : K \otimes K \to k$ such that $\sigma = \tau_\pi$.
As noted in \cite{bb}, we have an injective algebra map
\begin{align*}
\theta : H^{\tau_\pi} & \longrightarrow {_\tau \! K}  \otimes K_{\tau^{-1}} \otimes H\\
[x] & \longmapsto \{\pi(x_1)\}  \otimes \langle \pi(x_3)\rangle \otimes x_{2} 
\end{align*}
Consider the linear map $\psi =\varepsilon \otimes \varepsilon : {_\tau \! K}  \otimes K_{\tau^{-1}} \longrightarrow k$. We have $(\psi \otimes {\rm id}) \circ \theta = {\rm id}$, and hence
is a coalgebra map. Let $\varphi_0 : H \longrightarrow B$ be an inner faithful
finite-dimensional  representation, and let $\varphi : H \longrightarrow B \times B^{\rm op}$, 
$x \longmapsto (\varphi_0(x), \varphi_0(S(x)))$. It is clear that
$\varphi$ is still inner faithful, and the second condition in the previous theorem 
is satisfied since $S^2$ is inner. Hence the previous
theorem ensures that  the representation $({\rm id} \otimes \varphi) \circ \theta : H \longrightarrow 
{_\tau \! K}  \otimes K_{\tau^{-1}}\otimes (B\times B^{\rm op})$ is inner faithful,
and we are done since ${_\tau \! K}  \otimes K_{\tau^{-1}}$ is finite-dimensional.
\end{proof}

To illustrate the previous theorem,  
let us  examine the case of multiparametric deformations of ${\rm GL}_n$ of \cite{ast}.
We assume that the base field is $\C$ until the end of the section.
Let $\mathbf p = (p_{ij})\in M_n(\C)$ be a multiplicatively antisymmetric matrix: 
$p_{ii}=1$ and $p_{ij}p_{ji}=1$ for all $i,j$.
The algebra $\mathcal O_{\mathbf p}({\rm GL}_n(\C))$ is the algebra presented by generators
$x_{ij}$, $y_{ij}$, $ 1\leq i,j\leq n$ submitted to the following relations ($1 \leq i,j,k,l \leq n$):
$$x_{kl}x_{ij} = p_{ki}p_{jl}x_{ij}x_{kl}, \quad y_{kl}y_{ij} = p_{ki}p_{jl}y_{ij}y_{kl}, \quad
y_{kl}x_{ij} = p_{ik}p_{lj}x_{ij}y_{kl}$$
$$\sum_{k=1}^n x_{ik}y_{jk} = \delta_{ij} = \sum_{k=1}^n x_{ki}y_{kj}$$ 
The presentation we have given avoids the use of the quantum determinant.
The algebra $ \mathcal O_{\mathbf p}({\rm GL}_n(\C))$ has a standard Hopf algebra structure, described
as follows:
$$\Delta(x_{ij})=\sum_k x_{ik} \otimes x_{kj}, \ \Delta(y_{ij})=\sum_k y_{ik} \otimes y_{kj}, \
\varepsilon(x_{ij}) =\delta_{ij} = \varepsilon(y_{ij}), \ S(x_{ij}) = y_{ji}, \ S(y_{ij})= x_{ji}$$
When $p_{ij}=1$ for all $i,j$, one gets the usual Hopf algebra
 $\mathcal O_{\mathbf p}({\rm GL}_n(\C))$.
 It is known that $\mathcal O_{\mathbf p}({\rm GL}_n(\C))$ is a cotwist of
$\mathcal O({\rm GL}_n(\C))$, by  a cotwist induced 
by the quotient Hopf algebra $\mathcal O((\C^\times)^n) \simeq \C[\mathbb Z^n]$ (see e.g. \cite{ast}).   

\begin{corollary}
 Let $\mathbf p = (p_{ij})\in M_n(\C)$ be a multiplicatively antisymmetric matrix whose entries are roots of unity. Then the Hopf algebra $ \mathcal O_{\mathbf p}({\rm GL}_n(\C))$ is inner linear.
\end{corollary}

\begin{proof}
Let $M$ be the (finite and cyclic) subgroup of $\mathbb C^\times$ generated by the elements $p_{ij}$.
The cotwist
 defining $\mathcal O_{\mathbf p}({\rm GL}_n(\C))$ is induced
 by the quotient Hopf algebra $\mathcal O(M^n)$, and thus Theorem \ref{cotwist}
 gives the result.
\end{proof}

When the multiplicatively antisymmetric matrix $\mathbf p$ satisfies
some appropriate additionnal conditions, we have similar deformations for various algebraic subgroups 
of ${\rm GL}_n$ and Theorem 6.6 holds true with similar proof. For example if  $p_{ij}=-1,$ $\forall i \not =j$, one gets a  Hopf algebra deformation of $\mathcal O({\rm O}(n, \C))$, 
found in Section 4 of \cite{bbc} and corresponding to
the quantum symmetry group of the hypercube, which therefore is inner linear.

There are however interesting situations for which Theorem \ref{cotwist}
is not sufficient to ensure inner linearity, and where
we need the full strength of Theorem 6.4, for example for the Hopf algebra
$\mathcal O_{-1}(SL_2(\C))$ (which is well-known not to be a cotwist of 
$\mathcal O(SL_2(\C))$).

\begin{corollary}
 The Hopf algebra $\mathcal O_{-1}(SL_n(\C))$ is inner linear for any $n\geq 1$.
\end{corollary}

\begin{proof}
 The case $n=1$ is trivial while if $n \geq 3$, similarly to the previous cases, $\mathcal O_{-1}(SL_n(\C))$
 is a cotwist of $\mathcal O(SL_n(\C))$ induced by a finite-dimensional quotient.
 So we assume that $n=2$. Similarly to \cite{z}, we a an algebra embedding
 $$\theta : \mathcal O_{-1}(SL_2(\C)) \hookrightarrow M_2(\C) \otimes \mathcal O(SL_2(\C))$$
 $$a \mapsto \sigma_1\otimes a, \ b \mapsto \sigma_2\otimes b, \
 c \mapsto \sigma_2\otimes c, \ d \mapsto \sigma_1\otimes d$$
 where 
 $$\sigma_1= \begin{pmatrix}
              0 & 1 \\
	      1 & 0
             \end{pmatrix}, \quad \sigma_2 = \begin{pmatrix}
	    i & 0 \\
	    0 & -i \end{pmatrix} 
$$ and $a,b,c,d$ denote the standard generators of both $\mathcal O(SL_2(\C))$ and 
$\mathcal O_{-1}(SL_2(\C))$. Let $\psi : M_2(\C) \rightarrow \C$ be the unique
linear map such that $\psi(1)=1=\psi(\sigma_1)=\psi(\sigma_2)=\psi(\sigma_1\sigma_2)$.
It is clear that  $(\psi \otimes {\rm id}) \circ \theta$
is a vector space isomorphism
$\mathcal O_{-1}(SL_2(\C)) \longrightarrow  \mathcal O(SL_2(\C))$, and it is not difficult to check that it is a coalgebra map.
Now choose topological generators $x,y$ of $SL_2(\C)$ and consider the inner faithful representation
(Example 5.2) 
$$\varphi : \mathcal O(SL_2(\C)) \rightarrow \C^4, \ f \mapsto (f(x), f(y), f(x^{-1}, f(y^{-1}))$$
The second condition in Theorem 6.4 is satisfied, and hence 
we get an inner faithful representation
$\mathcal O_{-1}(SL_2(\C)) \rightarrow M_2(\C) \otimes \C^4$.
\end{proof}

\section{Tensor and free product of representations}

The Hopf image does not behave  well with respect to the tensor or free product:
the Hopf image of a tensor or free product is not necessarily 
the tensor or free product of the Hopf images. Here is what can be said in general.

\begin{proposition}
Let $H$ and $L$ be  Hopf algebras and
 let  $\pi : H \longrightarrow A$ and $\varphi : L \longrightarrow B$
 be algebra maps. Then there are surjective Hopf algebra maps
 $$H_\pi \otimes L_\varphi \longrightarrow (H \otimes L)_{\pi \otimes \varphi}
 \quad {\rm and} \quad H_\pi * L_\varphi \longrightarrow (H * L)_{\pi * \varphi}$$ 
\end{proposition}

\begin{proof}
 We consider the Hopf algebra factorizations
$$\xymatrix{
H\otimes L \ar[rr]^{\pi \otimes \varphi} \ar[dr]_{p_\pi \otimes p_\varphi} & & A \otimes B  \\
& H_\pi \otimes L_\varphi \ar[ur]_{\tilde{\pi} \otimes \tilde{\varphi}} &
} \quad 
\xymatrix{
H* L \ar[rr]^{\pi * \varphi} \ar[dr]_{p_\pi * p_\varphi} & & A * B  \\
& H_\pi * L_\varphi \ar[ur]_{\tilde{\pi} * \tilde{\varphi}} &
}
$$
and the  universal property of the Hopf image yields the announced Hopf algebra maps, which are clearly surjective.
\end{proof}

The morphisms in the proposition are not injective in general, as the following example shows.

\begin{example}{\rm 
Let $\z_n$ be the cyclic group of order $n$ and let $x \in \z_n$ be a generator.
Let $\pi : \C[\z_n] \longrightarrow \C$ be the unique algebra map with
$\pi(x) = \omega$, where $\omega$ is a primitive root of unity of order $n$.
The representation $\pi$ is inner faithful because its restriction to $\z_n$ is faithful, 
hence $\C[\z_n]_\pi = \C[\z_n]$. However the representations
$\pi \otimes \pi$ and $\pi*\pi$ are not inner faithful because 
the groups $\z_n \times \z_n$ and $\z_n*\z_n$ do not inject in $\C^*$.} 
\end{example}

For the tensor product, a faithfulness assumption on one 
of the algebra morphisms enables one to describe the Hopf image easily. We begin
with a lemma. We use the notation of Section 2.

\begin{lemma}
 Let $H$ and $L$ be Hopf algebras and
 let  $\pi : H \longrightarrow A$ and $\varphi : L \longrightarrow B$
 be algebra maps. Then
$$ I_{\pi \otimes \varphi}  = \bigcap_{g \in F} {\rm Ker}(\pi^g \otimes \varphi^g)$$
\end{lemma}

\begin{proof}
We first show that for any
$g \in F$, there exists a linear isomorphism
$T_g : A^g \otimes B^g \longrightarrow (A\otimes B)^g$
such that $(\pi \otimes \varphi)^g = T_g \circ (\pi^g \otimes \varphi^g)$.
We prove this by induction on $l(g)$. This is clear if $l(g)\leq 1$, so we assume that
$g=hh'$ with $l(g)>l(h)\geq 1$ and $l(g) >l(h')\geq 1$. 
Let $x \in H$ and $y \in L$. We have, using the induction assumption,
\begin{align*}
(\pi \otimes \varphi)^{hh'}(x \otimes y) & =
(\pi \otimes \varphi)^{h}(x_1 \otimes y_1) \otimes (\pi \otimes \varphi)^{h'}(x_2 \otimes y_2) \\
& = (T_h \otimes T_{h'}) (\pi^h(x_1) \otimes \varphi^h(y_1) \otimes \pi^{h'}(x_2) \otimes \varphi^{h'}(y_2))\\
& = T_{hh'}\circ (\pi^{hh'} \otimes \varphi^{hh'})( x \otimes y) 
\end{align*} 
and hence we have the desired result.
Thus for $g \in F$, we have ${\rm Ker}((\pi\otimes \varphi)^g) = {\rm Ker}(\pi^g \otimes \varphi^g)$,
and we have our result.
\end{proof}

\begin{proposition}
 Let $H$ and $L$ be Hopf algebras and
 let  $\pi : H \longrightarrow A$ and $\varphi : L \longrightarrow B$
 be algebra maps. Assume that  $\pi$ is  faithful and that the antipode of $H$ is injective. 
 Then
 $(H \otimes L)_{\pi \otimes \varphi}\cong H \otimes L_\varphi $. 
 In particular if $\pi$ is faithful and $\varphi$ is inner faithful, then
 $\pi \otimes \varphi$ is inner faithful.
\end{proposition}

\begin{proof}
 It is sufficent to show that
$I_{\pi \otimes \varphi} = H \otimes I_\varphi$, since 
$(H \otimes L) / (H \otimes I_\varphi) \simeq H \otimes (L/I_\varphi)$.
For any $g \in F$, the representation $\pi^g$ is faithful since $\Delta$ and $S$ are injective, so
${\rm Ker}(\pi^g \otimes \varphi^g) = H \otimes {\rm Ker}(\varphi^g)$. By the previous lemma we have 
\begin{align*}
 I_{\pi \otimes \varphi} & = \bigcap_{g \in F} {\rm Ker}(\pi^g \otimes \varphi^g) \\
& =  \bigcap_{g \in F} H \otimes {\rm Ker}(\varphi^g) = 
H \otimes(\bigcap_{g \in F} {\rm Ker}(\varphi^g))\\
&  = H \otimes I_\varphi
\end{align*}
and we have our result. The last assertion is immediate.
\end{proof}

\begin{remark}
 {\rm Of course the assumption of the injectivity of the antipode is satisfied 
 in the most interesting cases. See \cite{sc00} for an example
 of a Hopf algebra having a non-injective antipode.}
\end{remark}

The result in the previous proposition naturally leads to the following question.

\begin{question}
 Let  $\pi : H \longrightarrow A$ be a representation of a Hopf algebra $H$ on 
 an algebra $A$. Does there exist a condition on $\pi$, weaker than faithfulness, that
 ensures that for any inner faithful representation $\varphi : L \longrightarrow B$ of a Hopf algebra $L$ on an algebra $B$, then the representation $\pi \otimes \varphi$ is inner faithful?
\end{question}

If a representation $\pi$ satisfies to the hypothetic condition 
of Question 7.6, then in particular  
$\pi \otimes \pi$ will be inner faithful. It seems to be simpler 
to study first representations with this weaker property, and this leads to the following definition.

\begin{definition}
 Let  $\pi : H \longrightarrow A$ be a representation of a Hopf algebra $H$ on 
 an algebra $A$. We say that $\pi$ is \textbf{projectively inner faithful} if 
 $\pi \otimes \pi$ is inner faithful.
\end{definition}

Once again the terminology is motivated by the discrete group case: a representation
$\pi : k[\Gamma] \longrightarrow A$ is projectively inner faithful if and only if the associated
group morphism $\Gamma \longrightarrow A^\times /k^*$ is faithful.
When $A = {\rm End}(V)$, this means that $\Gamma$ embeds into the projective linear group ${\rm PGL}(V)$.

The example of function algebras shows that it is in general difficult
to decide when an inner faithful representation is projectively inner faithful. 
This specific example leads to the following definition.
This is a purely discrete group theoretic definition, but 
of course we have in mind discrete groups embedded as dense subgroups of algebraic groups.

\begin{definition}
 Let $\Gamma$ be a discrete group and let $\{s_1, \ldots , s_n\}$ be a family
 of generators of $\Gamma$. We say that $\{s_1, \ldots , s_n\}$ is a \textbf{family of projective
 generators} of $\Gamma$ if the group $\Gamma \times \Gamma$ is generated
 by the elements $(s_i,s_j)$, $1 \leq i,j \leq n$.
\end{definition}

If $\Gamma = \langle s_1, \ldots , s_n \rangle$ is a dense subgroup of an algebraic group
$G$, then $\{s_1, \ldots , s_n\}$ is a family of projective generators of $\Gamma$ if and only
if the corresponding representation $\mathcal O(G) \longrightarrow \C^n$ of Example 5.2 is 
projectively inner faithful.

Of course $\{s_1, \ldots , s_n\}$ is a family of projective generators if $1 \in \{s_1, \ldots s_n\}$.
It would be interesting to have more examples and to characterize
the family of projective generators of a group.

\medskip

We end the section by noting that there is no analogous result to Proposition 7.4. 
in the free product case, as shown by the following example.

\begin{example}
 {\rm Let 
 $\pi : \C[\z_n] \longrightarrow \C$ be the inner faithful representation
 of Example 7.2. Then $\pi * {\rm id} : \C[\z_n]*\C[\z_n]  \longrightarrow \C[\z_n]$ is not inner faithful by Proposition 2.12, since $\C[\z_n]*\C[\z_n]\simeq \C[\z_n*\z_n]$ is a noncommutative algebra
while $\C[\z_n]$ is commutative. }
\end{example}

\section{Tannaka duality}
Tannaka duality studies the interplays between 
a Hopf algebra and its category of comodules. 
In this section we formulate some Tannaka type results 
for Hopf images. These results are used in \cite{bbs}, in the study of quantum permutation groups associated to complex Hadamard matrices.

Let $H$ be a Hopf algebra and let $U,V$ be $H$-comodules.
The coaction on $U$ is denoted $\alpha_U : U \longrightarrow U \otimes H$.
The set of $H$-comodule morphisms from $U$ to $V$ is denoted 
${\rm Hom}_H(U,V)$. If $f : H \longrightarrow L$ is a Hopf algebra map,
then $f$ induces natural $L$-comodule structures on $U$ and $V$
(the resulting comodules still being denoted $U$ and $V$ if no confusion arises), with
$${\rm Hom}_H(U,V) \subset {\rm Hom}_L(U,V)$$

We now introduce a key space for the computation 
of Hom spaces of a Hopf image.

\begin{definition}
Let $\pi : H \longrightarrow A$ be a representation
of a Hopf algebra $H$ on an algebra $A$, and let $U,V$ be some $H$-comodules.
The space of $\pi$-morphisms from $U$ to $V$ is defined by
$${\rm Hom}(U_\pi,V_\pi) = 
\{ f \in {\rm Hom}_k(U,V) \ | \ (1_V \otimes \pi) \circ \alpha_V \circ f = 
(1_V \otimes \pi)\circ (f \otimes 1_H)  \circ \alpha_U\}$$
\end{definition}

The idea is that the space of $\pi$-morphisms is a more concrete one than the
space of $H_\pi$-comodules (at least if the algebra $A$ is a more concrete one than
$H$), and hence should be easier to describe. In fact it contains all the necessary information.

\begin{theorem}
 Let $\pi : H \longrightarrow A$ be a representation
of a Hopf algebra $H$ on an algebra $A$. Let $U$ and $V$ be finite dimensional 
$H$-comodules and let
$(L,q, \varphi)$ be a factorization of $\pi$. Then we have
$${\rm Hom}_H(U,V) \subset {\rm Hom}_L(U,V) \subset
{\rm Hom}_{H_\pi}(U,V)= {\rm Hom}(U_\pi,V_\pi)$$
\end{theorem}

\begin{proof}
 The first three inclusions on the left arise from the Hopf algebra maps $H \rightarrow L \rightarrow H_\pi$. Let $f \in {\rm Hom}_{H_\pi}(U,V)$.
 Then we have, with the notations of Theorem 2.1: 
 \begin{align*}
 & (f \otimes 1_{H_\pi}) \circ (1_U \otimes p) \circ \alpha_U = (1_V \otimes p)
 \circ \alpha_V \circ f \\
 \Rightarrow & (1_V \otimes p)\circ (f \otimes 1_H) \circ \alpha_U = (1_U \otimes p)
 \circ \alpha_V \circ f \\
 \Rightarrow &(1_V \otimes \pi) \circ \alpha_V \circ f = 
(1_V \otimes \pi)\circ (f \otimes 1_H)  \circ \alpha_U
 \end{align*}
 Hence $f \in {\rm Hom}(U_\pi, V_\pi)$ and ${\rm Hom}_{H_\pi}(U,V) \subset {\rm Hom}(U_\pi,V_\pi)$. 
 Assume conversely that $f \in {\rm Hom}(U_\pi, V_\pi)$, and let
 $e_1, \ldots , e_m$ and $e'_1, \ldots , e'_m$ be respective bases of $U$ and $V$ with
 $$\alpha_U(e_i) = \sum_{k=1}^n e_k \otimes u_{ki} \quad {\rm and } \quad
 \alpha_V(e'_j) = \sum_{l=1}^m e'_l \otimes v_{lj}$$
 Let $(\lambda_{ij}) \in M_{m,n}(k)$ be such that $f(e_i) = \sum_j \lambda_{ji}e'_j$.
 Then since $f \in {\rm Hom}(U_\pi, V_\pi)$, we have for $1 \leq i\leq n$, $1 \leq l\leq m$: 
 $$\sum_{j=1}^m \pi(v_{lj})\lambda_{ji} = \sum_{k=1}^n \lambda_{lk} \pi(u_{ki})$$
 and hence
 $$P_{li} := \sum_{j=1}^m v_{lj}\lambda_{ji} - \sum_{k=1}^n \lambda_{lk} u_{ki} \in {\rm Ker}(\pi)$$
 Let $I$ be the ideal generated by the elements $P_{li}$.
 We have 
 $$\Delta(P_{li}) = \sum_{j=1}^m v_{lj} \otimes P_{ji} + \sum_{k=1}^n P_{lk} \otimes u_{ki}$$ 
 and hence
 that $I$ is a bi-ideal. Multiplying $P_{il}$ on the left by $S(v_{rl})$ and on the right by
 $S(u_{is})$ and summing over $i$ and $l$ gives
 $$\sum_{i}\lambda_{ri} S(u_{is})- \sum_{l=1}^n S(v_{rl})\lambda_{ls} = S(-P_{rs}) \in I$$
Thus $I$ is a Hopf ideal and hence $I \subset I_\pi = {\rm Ker}(p)$. Thus  
we have, for $1 \leq i\leq n$, $1 \leq l\leq m$: 
 $$\sum_{j=1}^m p(v_{lj})\lambda_{ji} = \sum_{k=1}^n \lambda_{lk} p(u_{ki})$$
 This  exactly means that $f \in {\rm Hom}_{H_\pi}(U,V)$, and we are done.
\end{proof}

\begin{corollary}
 Let $\pi : H \longrightarrow A$ be a representation
of a Hopf algebra $H$ on an algebra $A$. If 
$\pi$ is inner faithful, then we have,
$${\rm Hom}_{H}(U,V)= {\rm Hom}(U_\pi,V_\pi)$$
 for any $H$-comodules  $U$ and $V$.
\end{corollary}

The converse of the corollary is not true in general. To see this, assume 
that $k$ has characteristic zero and consider 
$H = \mathcal O({\rm SL}_2(k))$ and $L = \mathcal O(B)$, with $B$ being
the Borel subgroup of ${\rm SL}_2(k)$ consisting of triangular matices.
The restriction map $\mathcal O({\rm SL}_2(k)) \longrightarrow \mathcal O(B)$
is not inner faithful because it is not faithful, and has
$\mathcal O(B)$ as Hopf image. One 
easily sees that for the irreducible representations of ${\rm SL}_2(k)$
(the symmetric powers of the fundamental representation), one 
has ${\rm Hom}_{{\rm SL}_2(k)}(U,V)= {\rm Hom}_B(U,V)$, and hence
by the cosemisimplicity of $\mathcal O({\rm SL}_2(k))$, this is true
for any representation of ${\rm SL}_2(k)$.

However we have a converse to the corollary if
both $H$ and the Hopf image are cosemisimple.

\begin{theorem}
 Let $\pi : H \longrightarrow A$ be a representation
of a cosemisimple Hopf algebra $H$ on an algebra $A$. Assume that the Hopf image 
$H_\pi$ is cosemisimple.
Then 
$\pi$ is inner faithful if and only if 
$${\rm Hom}_{H}(U,V)= {\rm Hom}(U_\pi,V_\pi)$$
for any simple $H$-comodules  $U$ and $V$.
\end{theorem}

\begin{proof}
The $\Rightarrow$ part is the previous corollary.
Conversely,  let us assume that
$${\rm Hom}_{H}(U,V)= {\rm Hom}(U_\pi,V_\pi)={\rm Hom}_{H_\pi}(U,V)$$
for any simple $H$-comodules  $U$ and $V$. 
Then the canonical projection $p: H \longrightarrow H_\pi$
induces an injection from the set of simple $H$-comodules
to the set of simple $H_\pi$-comodules, and using the respective Peter-Weyl decompositions of $H$ and $H_\pi$, we see that $p$ is injective, and hence is an isomorphism. 
\end{proof}

It seems to be difficult in general
to ensure that the Hopf image is cosemisimple. 
However this is automatically true when one works 
in the category of compact Hopf algebras (in the sense of Section 2),
and hence we get the following interesting characterization of inner faithfulness
for $*$-representations.

\begin{theorem}
 Let $\pi : H \longrightarrow A$ be a $*$-representation
of a compact Hopf algebra $H$ on a $*$-algebra $A$. 
Then 
$\pi$ is inner faithful if and only if 
$${\rm Hom}_{H}(U,V)= {\rm Hom}(U_\pi,V_\pi)$$
for any simple $H$-comodules  $U$ and $V$.
\end{theorem}

\begin{proof}
 The proof is done by the straightforward adaptation of the arguments of Theorem 8.2 and Theorem 8.4
 to the $*$-case.
\end{proof}

There is also a useful variant of the above theorem, when 
$H$ is finitely generated. 
Let $U$ be a finite-dimensional $H$-comodule. To any $x \in \N * \N = \langle \alpha, \beta\rangle $, we 
associate an $H$-comodule $U^x$ as follows: 
$U^1= \C$, $U^\alpha = U$, $U^\beta = U^*$, $U^{\alpha \beta} = U \otimes U^*$, 
$U^{\beta \alpha} = U^* \otimes U$ and so on.
We say that the $H$-comodule $U$ is \textbf{faithful} if any finite
finite-dimensional comodule is a subquotient of a direct sum of objects
of type $U^x$. Assuming that $H$ is compact, this is equivalent to saying
that any simple $H$-comodule is a subobject of some $U^x$. 
If $u=(u_{ij})$ is the matrix of coefficients of $U$, it is known
that $U$ is faithful if and only if the coefficients $u_{ij}$ and $S(u_{ij})$ generate $H$ as an algebra.

\begin{theorem}
 Let $\pi : H \longrightarrow A$ be a $*$-representation
of a compact Hopf algebra $H$ on a $*$-algebra $A$, and let 
$U$ be a faithful $H$-comodule.
Then $\pi$ is inner faithful if and only if
 $${\rm Hom}_{H}(\C,U^x)= {\rm Hom}(\C_\pi,U^x_\pi)$$
 for any $x \in \N * \N$. 
\end{theorem}

\begin{proof}
 The $\Rightarrow$ part is again  Corollary 8.3.
Conversely, assume that ${\rm Hom}_{H}(\C,U^x)= {\rm Hom}(\C_\pi,U^x_\pi)$
 for any $x \in \N * \N$. By duality theory in monoidal categories, we get
 $${\rm Hom}_{H}(U^x,U^y)= {\rm Hom}(U^x_\pi,U^y_\pi) \ (={\rm Hom}_{H_\pi}(U^x,U^y))$$
for any $x,y \in \N*\N$. Hence by Lemma 5.3 in \cite{ba99}
we see that the canonical morphism $p:H \longrightarrow H_\pi$ is an isomorphism,
which shows that $\pi$ is inner faithful.
\end{proof}

When the faithful $H$-comodule is self-dual
(for example in the case of quantum permutation algebras as in \cite{bbs}), 
the previous theorem has the following simpler form. 

\begin{theorem}
 Let $\pi : H \longrightarrow A$ be a $*$-representation
of a compact Hopf algebra $H$ on a $*$-algebra $A$, and let 
$U$ be a faithful self-dual $H$-comodule.
Then $\pi$ is inner faithful if and only if
 $${\rm Hom}_{H}(\C,U^{\otimes n})= {\rm Hom}(\C_\pi,U^{\otimes n}_\pi)$$
 for any $n \in \N$. 
\end{theorem}

\begin{proof}
 All the comodules $U^x$ of the previous theorem are isomorphic
 with $U^{\otimes n}$ for some $n$, and hence the result follows.
\end{proof}

\section{Hopf algebras with small corepresentation level}

This section gives a concrete application of the inner faithfulness criterion
of the previous section to compact Hopf algebras having all their
simple comodules of dimension smaller than $2$.
The Hopf algebras we consider arise in the study of
$4 \times 4$ Hadmard matrices (\cite{bn,bbs}).  
 Let us begin 
with some vocabulary.

\begin{definition}
 Let $H$ be a Hopf algebra and let $n \in \mathbb N^*$. We say that 
 $H$ has \textbf{corepresentation level $n$} if there exists
 a simple $H$-comodule $V$ with $\dim(V)=n$, and if any simple
 $H$-comodule has dimension smaller than $n$.
 In this case we put $cl(H) = n$. If the set 
 of dimensions of simple $H$-comodules is not bounded, we put $cl(H) = \infty$. 
\end{definition}

Pointed Hopf algebras are exactly the Hopf algebras
with $cl(H)=1$.
The finite groups $\Gamma$ having all their irreducible representations
of dimension $\leq 2$ (in characteristic zero) are described in \cite{am}, corresponding
to the function Hopf algebras $k^\Gamma$ such that $cl(k^\Gamma) \leq 2$.

We begin with a general result (Theorem 9.2) that ensures that a representation of 
a compact Hopf algebra with corepresentation level $2$ is inner faithful.
Then this result will be used to construct inner faithful representations
of some concrete Hopf algebras in Theorems 9.3 and 9.4.

\begin{theorem}
Let $H$ be a compact Hopf algebra with $cl(H)=2$.
Let $\Gamma$ be the group of group-like elements of $H$ and let
$\Lambda$ be the set of isomorphism classes of simple two-dimensional
$H$-comodules. For each $\lambda \in \Lambda$, fix a matrix $u^\lambda =(u_{ij}^\lambda)\in M_2(H)$
of corresponding coefficients.
Let $\pi : H \longrightarrow A$ be a $*$-representation. Assume that the following
conditions are fulfilled.
\begin{enumerate}
 \item $\pi_{|\Gamma}$ is injective.
 \item $\forall \lambda \in \Lambda$, $\pi(u_{11}^\lambda)= 0=\pi(u_{22}^\lambda)$.
 \item $\forall \lambda \in \Lambda$, $\forall g \in \Gamma$,
 $\pi(u_{12}^\lambda)$ and $\pi(g)$ are linearly independent.
 \item $\forall \lambda, \mu \in \Lambda$, $\pi(u_{12}^\lambda)$ and $\pi(u_{21}^\mu)$
 are linearly independent.
 \item $\forall \lambda, \mu \in \Lambda$ with $\lambda \not = \mu$, $\pi(u_{12}^\lambda)$ and $\pi(u_{12}^\mu)$
 are linearly independent.
\end{enumerate}
Then $\pi$ is inner faithful. 
\end{theorem}

\begin{proof}
We have to show, by Theorem 8.4, that
 ${\rm Hom}_{H}(U,V)= {\rm Hom}(U_\pi,V_\pi)$
for any simple $H$-comodules  $U$ and $V$.
If $U$ and $V$ have dimension $1$, this is ensured by condition (1).
Assume that $U$ has dimension $2$, corresponding to $\lambda \in \Lambda$ and that $V$ has dimension $1$, corresponding to $g \in G$. Then  ${\rm Hom}(U_\pi,V_\pi)$ is identified
with the spaces of matrices $t=(t_1,t_2) \in M_{1,2}(\C)$ such that
$$(t_1,t_2)
\begin{pmatrix} \pi(u_{11}^\lambda) & \pi(u_{12}^\lambda) \\
\pi(u_{21}^\lambda) & \pi(u_{22}^\lambda) \end{pmatrix}
=(t_1\pi(g),t_2\pi(g))$$
Using assumptions (2) and (3), we see that ${\rm Hom}(U_\pi,V_\pi) = \{0\} = {\rm Hom}_H(U,V)$.
Using a similar argument, we see that the same conclusion holds if
$\dim(U)=1$ and $\dim(V)=2$. Assume now that $\dim(U)=2=\dim(V)$, corresponding respectively to 
$\lambda, \mu \in \Lambda$. Then ${\rm Hom}(U_\pi,V_\pi)$ is identified
with the space of matrices $t=(t_{ij}) \in M_{2}(\C)$ such that
$$\begin{pmatrix} \pi(u_{11}^\mu) & \pi(u_{12}^\mu) \\
\pi(u_{21}^\mu) & \pi(u_{22}^\mu) \end{pmatrix} 
\begin{pmatrix} t_{11} & t_{12} \\
t_{21} & t_{22} \end{pmatrix} =
\begin{pmatrix} t_{11} & t_{12} \\
t_{21} & t_{22} \end{pmatrix}
\begin{pmatrix} \pi(u_{11}^\lambda) & \pi(u_{12}^\lambda) \\
\pi(u_{21}^\lambda) & \pi(u_{22}^\lambda) \end{pmatrix}$$
Using assumptions (2) and (4), we see that $t_{21}=t_{12}=0$ and
if $\lambda = \mu$ we have $t_{11}=t_{22}$, and hence
${\rm Hom}(U_\pi,V_\pi)  = {\rm Hom}_H(U,V) = \C$.
If $\lambda \not = \mu$, we have $t_{11}=t_{22}=0$ by (5)
and hence
${\rm Hom}(U_\pi,V_\pi) = (0) = {\rm Hom}_H(U,V)$: this concludes the proof.
\end{proof}

We now use Theorem 9.2 to provide inner faithful representations
of a class of Hopf algebras considered in \cite{bb}, Section 7.
We refer the reader to \cite{bb} for the precise origins of these Hopf algebras, which we present now.
First we have the Hopf algebra $A_h(2)$: 
this is the universal $*$-algebra
 presented by generators
$(v_{ij})_{1\leq i,j\leq 2}$ and relations:
\begin{enumerate}
 \item The matrix  $v=(v_{ij})$ is orthogonal (with $v_{ij}^*=v_{ij}$).
 \item $v_{ij}v_{ik}= v_{ik}v_{kj} = 0=v_{ji}v_{ki} = v_{ki}v_{ji}$ if $j \not = k$. 
\end{enumerate}
The standard formulae 
$$\Delta(v_{ij}) = \sum_kv_{ik}\otimes v_{kj}, \quad \varepsilon(v_{ij})=\delta_{ij}, \
S(v_{ij})=v_{ji}$$
endow $A_h(2)$ with a compact Hopf algebra structure.
We have $cl(A_h(2))=2$, and the fusion rules of the corepresentations
of $A_h(2)$ are those of the orthogonal group O$(2)$. 

The Hopf algebra $A_h(2)$ has a series of finite-dimensional quotients
defined as follows.

For $m \in \mathbb N^*$ and $e= \pm 1$, the compact Hopf algebra $A(2m,{e})$ is the quotient
of $A_h(2)$ by the relations
$$(v_{11}v_{22})^m = (v_{22}v_{11})^m, \ 
(v_{12}v_{21})^m = e(v_{21}v_{12})^m$$
Similarly, for $m \in \mathbb N$ and $e=\pm 1$,
the compact Hopf algebra $A(2m+1,{e})$ is  the quotient
of $A_h(2)$ by the relations
$$(v_{11}v_{22})^mv_{11} = (v_{22}v_{11})^mv_{22}, \ 
(v_{12}v_{21})^m v_{12}= e(v_{21}v_{12})^mv_{21}$$
For any $k \in \mathbb N^*$, we have $\dim(A(k,e))=4k$
and $cl(A(k,e))=2$ if $k\geq 2$.

\begin{theorem}
 Let $q \in \C^*$ with $q \not =1$ and $|q|=1$.
 Then there exists a unique $*$-representation
 $\pi_q : A_h(2) \longrightarrow M_2(\C)$ such that
 $$\pi_q(v_{11}) = 0=\pi_q(v_{22}), \ 
 \pi_q(v_{12})=
 \begin{pmatrix}
  0 & q^{-1} \\ q & 0
 \end{pmatrix}, \
\pi_q(v_{21}) = 
\begin{pmatrix}
 0 & 1 \\
 1 & 0
\end{pmatrix}$$
\begin{enumerate}
\item If $q$ is not a root of unity, then the $*$-representation $\pi_q$ is inner faithful.
\item If $q$ has order $2m+1$, then $A_h(2)_{\pi_q} \simeq A(2m+1,1)$.
\item If $q$ has order $4m$, then $A_h(2)_{\pi_q} \simeq A(2m,-1)$.
\item If $q$ has order $4m+2$, then $A_h(2)_{\pi_q} \simeq A(2m+1,-1)$.
\end{enumerate}
\end{theorem}

\begin{proof} It is straightforward to check the existence  of the $*$-representation $\pi_q$.
For $m \geq 1$, we have
$$\pi((v_{11}v_{22})^m)=0 =  \pi((v_{22}v_{11})^m),$$
$$\pi((v_{12}v_{21})^m) = \begin{pmatrix}
  q^{-m} & 0\\ 0 & q^m 
 \end{pmatrix}, \ \pi((v_{21}v_{12})^m) = \begin{pmatrix}
  q^{m} & 0\\ 0 & q^{-m} 
 \end{pmatrix}$$
 and for $m \geq 0$ we have
 $$\pi((v_{11}v_{22})^mv_{11})=0 =  \pi((v_{22}v_{11})^m)v_{22},$$
$$\pi((v_{12}v_{21})^mv_{12}) = \begin{pmatrix}
  0 & q^{-m-1}\\  q^{m+1} & 0 
 \end{pmatrix}, \ \pi((v_{21}v_{12})^mv_{21}) = \begin{pmatrix}
 0&  q^{m}\\ q^{-m} & 0 
 \end{pmatrix}$$
 
 (1) We have $cl(A_h(2)) = 2$ and the simple $A_h(2)$-comodules are as follows: there is only one non-trivial
 one dimensional comodule, corresponding to the group-like $d= v_{11}^2-v_{12}^2$, and 
 a family of simple 2-dimensional comodules $V_n$, $n\geq 1$, corresponding to 
 the simple subcoalgebras
 $$C(2m) = \C(v_{11}v_{22})^m \oplus \C(v_{12}v_{21})^m \oplus \C(v_{21}v_{12})^m
 \oplus \C(v_{22}v_{11})^m, \ m\geq 1, \ {\rm and} $$
   $$C(2m+1) = \C(v_{11}v_{22})^mv_{11} \oplus \C(v_{12}v_{21})^mv_{12} \oplus \C(v_{21}v_{12})^mv_{21}
 \oplus \C(v_{22}v_{11})^mv_{22}, \ m\geq 0$$
 It is straightforward to check that if $q$ is not a root of unity, the 
 representation $\pi_q$ satisfies to the conditions of Theorem 9.2, and hence we have our result.

(2) Assume that $q$ is a root of unity of order $2m+1$. Then the representation
$\pi_q$ induces a $*$-representation $\tilde{\pi_q} : A(2m+1,1) \longrightarrow M_2(\C)$, and we have to show that $\tilde{\pi_q}$ is inner faithful.
For notational simplicity, the projections of elements of $A_h(2)$ in $A(2m+1,1)$
are denoted by the same symbol.
We have $cl(A(2m+1,1)) \leq 2$ and the simple $A(2m+1,1)$-comodules are as follows: there are 3 non-trivial
 $1$-dimensional comodules, corresponding to the group-like elements
  $$d= v_{11}^2-v_{12}^2, \ g= (v_{11}v_{22})^mv_{11} + (v_{12}v_{21})^mv_{12}, \
  h= (v_{11}v_{22})^mv_{11} - (v_{12}v_{21})^mv_{12}$$ 
 a family of $2m$ simple 2-dimensional comodules $V_1, \ldots , V_{2m}$, corresponding to 
 the simple subcoalgebras $C(k)$ just defined above.
 It is then a direct verification to check that the 
 representation $\tilde{\pi_q}$ satisfies to the conditions of Theorem 9.2, and hence is inner faithful.
 
 (3) Assume that $q$ is a root of unity of order $4m$ ($m\geq 1$). Then the representation
$\pi_q$ induces a $*$-representation $\tilde{\pi_q} : A(2m,-1) \longrightarrow M_2(\C)$, and we have to show that $\tilde{\pi_q}$ is inner faithful. We have
$cl(A(2m,-1)) =2$ and the simple $A(2m,-1)$-comodules are as follows: there are 3 non-trivial
 $1$-dimensional comodules, corresponding to the group-like elements
  $$d= v_{11}^2-v_{12}^2, \ g= (v_{11}v_{22})^m + i(v_{12}v_{21})^m, \
  h= (v_{11}v_{22})^m - i(v_{12}v_{21})^m$$ 
 a family of  $(2m-1)$ simple 2-dimensional comodules $V_1, \ldots , V_{2m-1}$, corresponding to 
 the simple subcoalgebras $C(k)$ defined above.
 It is  a direct verification to check that the 
 representation $\tilde{\pi_q}$ satisfies to the conditions of Theorem 9.2, and hence is inner faithful.
 
 (4) Assume finally that $q$ is a root of unity of order $4m+2$ ($m\geq 1$). Then the representation
$\pi_q$ induces a $*$-representation $\tilde{\pi_q} : A(2m+1,-1) \longrightarrow M_2(\C)$, and we have to show that $\tilde{\pi_q}$ is inner faithful. We have
$cl(A(2m+1,-1)) \leq 2$ and the simple $A(2m+1,-1)$-comodules are as follows: there are 3 non-trivial
 $1$-dimensional comodules, corresponding to the group-like elements
  $$d= v_{11}^2-v_{12}^2, \ g= (v_{11}v_{22})^mv_{11} + i(v_{12}v_{21})^mv_{12}, \
  h= (v_{11}v_{22})^mv_{11} - i(v_{12}v_{21})^mv_{12}$$ 
 a family of  $2m$ simple 2-dimensional comodules $V_1, \ldots , V_{2m}$, corresponding to 
 the simple subcoalgebras $C(k)$ defined above.
 It is  again a direct verification to check that the 
 representation $\tilde{\pi_q}$ satisfies to the conditions of Theorem 9.2, and hence is inner faithful.
\end{proof}

The Hopf algebra $A_h(2)$ is in fact a particular case 
of a construction given in \cite{bi2}, Example 2.5, that we describe now, 
and the inner faithful representation of the previous theorem
has a natural generalization.

Let $\Gamma$ be a group and let $n \in \mathbb N^*$.
Let $A_n(\Gamma)$ be the algebra presented by  generators
$a_{ij}(g)$, $1 \leq i,j \leq n$, $g \in \Gamma$, and   
relations ($1 \leq i,j,k \leq n$ ;   $g,h \in \Gamma$):
$$a_{ij}(g) a_{ik}(h) = \delta_{jk} a_{ij}(gh) \quad ; \quad
a_{ji}(g) a_{ki}(h) = \delta_{jk} a_{ji}(gh) \quad ; \quad
\sum_{l=1}^n a_{il}(1) = 1 = \sum_{l=1}^n a_{li}(1).$$ 
Then $A_n(\Gamma)$ is a compact Hopf algebra, with:
$$ a_{ij}(g)^* = a_{ij}(g^{-1}) \ ; \ 
\Delta(a_{ij}(g)) = \sum_{k=1}^n a_{ik}(g) \otimes a_{kj}(g) \ ; \
\varepsilon(a_{ij}(g)) = \delta_{ij} \ ; \
S(a_{ij}(g)) = a_{ji}(g^{-1}).$$
The Hopf algebra $A_h(2)$ is $A_2(\z_2)$. 
 
We need some notation to state
a generalization of the first part of Theorem 9.3.
We consider the group free product $\Gamma * \Gamma$, with the canonical
morphisms still denoted $\nu_1,\nu_2 : \Gamma \longrightarrow \Gamma * \Gamma$.
The canonical involutive group automorphism of $\Gamma*\Gamma$ is denoted by 
$\tau$, with $\tau \circ \nu_1 = \nu_2$ and $\tau \circ \nu_2 = \nu_1$.
%An element $x \in  G*G$ is the unit element if and only if
%$\tau(x) = x$.

\begin{theorem}
 Let $\Gamma$ be a group, let $A$ be a $*$-algebra, and 
 let $\pi_0 : \Gamma * \Gamma \longrightarrow A^\times$ be a group morphism
 into the group of unitary elements of $A$, such that 
  $\forall x \in \Gamma * \Gamma \setminus \{1\}$, we have $\pi_0(x) \not \in \C1$.  
Then $\pi_0$ induces an inner faithful $*$-representation
$\pi : A_2(\Gamma) \longrightarrow A$ such that
$$\forall g \in \Gamma, \ \pi(a_{11}(g)) = 0 = \pi(a_{22}(g)), \ \pi(a_{12}(g)) = \pi_0(\nu_1(g)), \ 
\pi(a_{21}(g)) = \pi_0(\nu_2(g))$$ 
In particular, if there exists a group embedding $\Gamma * \Gamma \subset {\rm PU}(n)$
for some $n \geq 1$, then the Hopf algebra $A_2(\Gamma)$ is inner linear.
\end{theorem}

\begin{proof}
 It is a direct verification to check that the above formulae define
 a $*$-representation $\pi : A_2(\Gamma) \longrightarrow A$. Recall from \cite{bi2}
 that we have an algebra isomorphism
 \begin{align*}
A_2(\Gamma) &\longrightarrow \C[(\Gamma *\Gamma) \times \z_2]\\
a_{ij}(g) &\longmapsto \nu_i(g)x_{ij} 
\end{align*}
where if $a$ is the generator of $\z_2$, $x_{11}= \frac{1+a}{2}=x_{22}$ and $x_{12} = \frac{1-a}{2} = x_{21}$ (of course we simply write $a$ for the element $(1,a)$ of   $(\Gamma *\Gamma) \times \z_2$,
and so on). We freely use this algebra identification in what follows.
Let us now recall the corepresentation theory
of $A_2(\Gamma)$ (Proposition 2.6 in \cite{bi2}).
We assume of course that the group $\Gamma$ is non-trivial
(otherwise the statement in the Theorem is trivial).
The only non trivial group-like element in 
$A_2(\Gamma)$ is $d=a_{11}(1)-a_{12}(1)$, for which $\pi(d)=-1$.
To any $x \in \Gamma *\Gamma \setminus\{1\}$ is associated a simple
comodule $U_x$ with corresponding matrix of coefficients defined by
$$A_{11}(x) = x x_{11} \ , \ A_{12}(x) = x x_{12} \ , \
A_{21}(x) = \tau(x)x_{21}\ , \  A_{22}(x) = \tau(x)x_{22}.$$ 
The comodules $U_x$ and $U_y$ are isomorphic if and only if
$x=y$ or $x=\tau(y)$. Any $2$-dimensional $A_2(\Gamma)$-comodule
is isomorphic to some $U_x$, and $cl(A_2(\Gamma))=2$. For 
$x \in \Gamma * \Gamma\setminus\{1\}$, we have
$$\pi(A_{11}(x)) = 0 = \pi( A_{22}(x)) \ , \ \pi(A_{12}(x)) = \pi_0(x) \ , \
\pi(A_{21}(x)) = \pi_0(\tau(x))$$ 
Then it is a straightforward verification to check 
that the conditions of Theorem 9.2 are fulfilled, using 
the properties of the group morphism $\pi_0$ (using also that $x=1$ if and only if $x=\tau(x)$), and we conclude that 
$\pi$ is inner faithful. 
\end{proof}

\end{document}